\documentclass[a4paper,11pt]{amsart}
\usepackage{hyperref}
\hypersetup{pdfborder=0 0 0}
\usepackage[utf8]{inputenc}
\usepackage{amssymb}
\usepackage{mathrsfs}
\usepackage[cal=boondox]{mathalfa}
\usepackage{tikz}
\usepackage{tikz-cd}
\usepackage{enumitem}
\usepackage{rotating}
\usepackage{MnSymbol}
\usepackage{relsize}
\usepackage{bbm}
\usepackage{tcolorbox}
\usepackage{xcolor}
\usepackage{mathtools}
\usepackage{musicography}
\usepackage[
backend=biber,
style=alphabetic,
sorting=anyt
]{biblatex}
\RequirePackage{doi}

\setlist{itemsep=1.0 pt}

\usepackage{textgreek}

\newcommand{\leftrarrows}{\mathrel{\raise.75ex\hbox{\oalign{%
  $\scriptstyle\leftarrow$\cr
  \vrule width0pt height.5ex$\hfil\scriptstyle\relbar$\cr}}}}
\newcommand{\lrightarrows}{\mathrel{\raise.75ex\hbox{\oalign{%
  $\scriptstyle\relbar$\hfil\cr
  $\scriptstyle\vrule width0pt height.5ex\smash\rightarrow$\cr}}}}
\newcommand{\Rrelbar}{\mathrel{\raise.75ex\hbox{\oalign{%
  $\scriptstyle\relbar$\cr
  \vrule width0pt height.5ex$\scriptstyle\relbar$}}}}

\makeatletter
\def\leftrightarrowsfill@{\arrowfill@\leftrarrows\Rrelbar\lrightarrows}
\newcommand{\xleftrightarrows}[2][]{\ext@arrow 3399\leftrightarrowsfill@{#1}{#2}}
\makeatother

\theoremstyle{plain}
\newtheorem{teo}{Theorem}[section]
\newtheorem{prop}[teo]{Proposition}
\newtheorem{lemma}[teo]{Lemma}
\newtheorem{cor}[teo]{Corollary}

\theoremstyle{definition}

\newtheorem{defi}[teo]{Definition}
\newtheorem{nota}[teo]{Notation}
\newtheorem{ass}[teo]{Assumption}
\newtheorem{set}[teo]{Setting}

\newtheorem{ex/}[teo]{Example}
\newenvironment{ex}
  {%
   \pushQED{\qed}\begin{ex/}}
  {\popQED\end{ex/}}

\newtheorem{rmk/}[teo]{Remark}
\newenvironment{rmk}
  {%
   \pushQED{\qed}\begin{rmk/}}
  {\popQED\end{rmk/}}

\numberwithin{equation}{section}
\newcommand{\Hom}{\mathrm{Hom}}

\newcommand{\Ker}{\mathrm{Ker}}
\newcommand{\End}{\mathrm{End}}

\newcommand{\Imm}{\mathrm{Im}}

\newcommand{\Q}{\mathbb{Q}}
\newcommand{\R}{\mathbb{R}}
\newcommand{\C}{\mathbb{C}}
\newcommand{\F}{\mathbb{F}}
\newcommand{\Z}{\mathbb{Z}}

\newcommand{\M}{\mathbb{M}}
\newcommand{\SSS}{\mathbb{S}}
\newcommand{\T}{\mathbb{T}}

\newcommand{\V}{\mathbb{V}}

\newcommand{\aid}{\mathfrak{a}}
\newcommand{\p}{\mathfrak{p}}
\newcommand{\q}{\mathfrak{q}}
\newcommand{\mmm}{\mathfrak{m}}
\newcommand{\cc}{\mathfrak{c}}
\newcommand{\lid}{\mathfrak{l}}
\newcommand{\fid}{\mathfrak{f}}
\newcommand{\Hf}{{\boldsymbol{f}}}
\newcommand{\Hg}{{\boldsymbol{g}}}
\newcommand{\Hh}{{\boldsymbol{h}}}
\newcommand{\HXi}{\pmb{\Xi}}
\newcommand{\Hxi}{\pmb{\xi}}
\newcommand{\Ad}{\mathbb{A}}
\newcommand{\BB}{\mathcal{B}}
\newcommand{\CC}{\mathcal{C}}
\newcommand{\DD}{\mathcal{D}}
\newcommand{\EE}{\mathcal{E}}

\newcommand{\HH}{\mathcal{H}}

\newcommand{\OO}{\mathcal{O}}

\newcommand{\RR}{\mathcal{R}}

\newcommand{\CS}{\mathcal{S}}
\newcommand{\W}{\mathcal{W}}

\newcommand{\UU}{\mathcal{U}}
\newcommand{\VV}{\mathcal{V}}

\newcommand{\SL}{\mathrm{SL}}
\newcommand{\GL}{\mathrm{GL}}

\newcommand{\re}{\mathrm{Re}}

\newcommand{\ord}{\mathrm{ord}}

\newcommand{\Gal}{\mathrm{Gal}}
\newcommand{\lcm}{\mathrm{lcm}}
\newcommand{\Rank}{\mathrm{Rank}}
\newcommand{\Sum}{\mathlarger{\sum}}
\newcommand{\Frac}{\mathrm{Frac}}
\newcommand{\Pic}{\mathrm{Pic}}
\newcommand{\Frob}{\mathrm{Frob}}
\newcommand{\Ind}{\mathrm{Ind}}

\newcommand{\Tr}{\mathrm{Tr}}
\newcommand{\Vol}{\mathrm{Vol}}

\newcommand{\CCC}{\hat{\CC}}
\newcommand{\Gcol}{{\boldsymbol{g}_{_{\mathrm{Col}}}}}
\newcommand{\Ghid}{{\boldsymbol{g}_{_{\mathrm{Hida}}}}}

\newcommand{\cyc}{\varepsilon_\mathrm{cyc}}

\newcommand{\Ls}{\mathscr{L}}

\newcommand{\ad}{\mathrm{Ad}}

\calclayout
\begin{document}

\title{Generalized triple product \texorpdfstring{$p$}{p}-adic \texorpdfstring{$L$}{L}-functions and rational points on elliptic curves}
\author{Luca Marannino}
\address{Sorbonne Université and Université Paris Cité, CNRS, IMJ-PRG, F-75005 Paris, France}
\email{marannino@imj-prg.fr}
\date{last update: \today}

\begin{abstract}
We generalize and simplify the constructions of \cite{DR2014} and \cite{Hsi2021} of an unbalanced triple product $p$-adic $L$-function $\mathscr{L}_p^f(\Hf,\Hg,\Hh)$ attached to a triple $(\Hf,\Hg,\Hh)$ of $p$-adic families of modular forms, allowing more flexibility for the choice of $\Hg$ and $\Hh$. 

Assuming that $\Hg$ and $\Hh$ are families of theta series of infinite $p$-slope, we prove a factorization of (an improvement of) $\mathscr{L}_p^f(\Hf,\Hg,\Hh)$ in terms of two anticyclotomic $p$-adic $L$-functions. As a corollary, when $\Hf$ specializes in weight $2$ to the newform attached to an elliptic curve $E$ over $\Q$ with multiplicative reduction at $p$, we relate certain Heegner points on $E$ to certain $p$-adic partial derivatives of $\mathscr{L}_p^f(\Hf,\Hg,\Hh)$ evaluated at the critical triple of weights $(2,1,1)$.
\end{abstract}

\maketitle

\tableofcontents

\section{Introduction and statement of the main results}
\subsection{The generalized unbalanced triple product \texorpdfstring{$p$}{p}-adic \texorpdfstring{$L$}{L}-function}
Let $p\geq 3$ be a rational prime. We fix an algebraic closure $\bar{\Q}$ of $\Q$, an algebraic closure $\bar{\Q}_p$ of $\Q_p$ together with an embedding $\iota_p:\bar{\Q}\hookrightarrow\bar{\Q}_p$ extending the canonical inclusion $\Q\hookrightarrow\Q_p$. All algebraic extensions of $\Q$ (resp. $\Q_p$) are viewed inside the corresponding fixed algebraic closures. We extend the $p$-adic absolute value $|\cdot|_p$ on $\Q_p$ (normalized so that $|p|_p=1/p$) to $\bar{\Q}_p$ in the unique possible way. We denote by $\C_p$ the completion of $\bar{\Q}_p$ with respect to this absolute value. It is well-known that $\C_p$ is itself algebraically closed.
We also fix an embedding $\iota_{\infty}:\bar{\Q}\hookrightarrow\C$ extending the canonical inclusion $\Q\hookrightarrow\C$ and we often omit the embeddings $\iota_p$ and $\iota_\infty$ from the notation. 

\medskip
Let $L/\Q_p$ be a finite extension and let $\Lambda:=\OO_L[\![1+p\Z_p]\!]$ be the corresponding Iwasawa algebra ($\OO_L$ being the ring of integers of $L$). Consider a new, $L$-rational, Hida family
\[
\Hf=\sum_{n=1}^{+\infty}a_n(\Hf)q^n\in \SSS^\ord(N_\Hf,\chi_\Hf,\Lambda)
\]
of tame level $N_\Hf$ ($p\nmid N_\Hf$) and tame character $\chi_\Hf$ of conductor dividing $N_\Hf$.

\medskip
Let also
\[
\Hg=\sum_{n=1}^{+\infty}a_n(\Hg)q^n\in\SSS_{\Omega_1}(M,\chi_\Hg,R_\Hg)\quad\text{and}\quad \Hh=\sum_{n=1}^{+\infty}a_n(\Hh)q^n\in\SSS_{\Omega_2}(M,\chi_\Hh,R_\Hh)
\]
be two generalized normalized $\Lambda$-adic eigenforms with $\chi_\Hf\cdot\chi_\Hg\cdot\chi_\Hh=\omega^{2a}$ for some integer $a$, where $\omega$ denotes the Teichmüller character modulo $p$ and $N_\Hf\mid M$. 

Our notion of generalized $\Lambda$-adic forms takes inspiration from \cite[Definition 2.16]{DR2014}. For a precise definition and for the explanation of the notation we refer to section \ref{generalized}. Here we just mention that we are not imposing any condition on $p$-slopes and that we are allowing the rings of coefficients $R_\Hg$ and $R_\Hh$ to be complete local noetherian flat $\Lambda$-algebras (not necessarily finite as $\Lambda$-algebras), having the same residue field as $\OO_L$.

\medskip
If $\Hg$ and $\Hh$ are Hida families, the works of Darmon-Rotger \cite{DR2014} and Hsieh \cite{Hsi2021} attach to the triple $(\Hf,\Hg,\Hh)$ a so-called $\Hf$-unbalanced square-root triple product $p$-adic $L$-function. It arises as an element 
\[
\mathscr{L}_p^f(\Hf,\Hg,\Hh)\in R_{\Hf\Hg\Hh}:=\Lambda\hat{\otimes}_{\OO_L}R_\Hg\hat{\otimes}_{\OO_L}R_\Hh\,,
\]
whose square interpolates the central values of the triple product $L$-functions attached to the specializations of $(\Hf,\Hg,\Hh)$ at $\Hf$-unbalanced triples of weights.

\medskip
More precisely, given two primitive Hida families $\Hg^\#$ and $\Hh^\#$ of respective tame level $N_\Hg$ and $N_\Hh$, Hsieh associates to the triple $(\Hf,\Hg^\#,\Hh^\#)$ a preferred choice of \emph{test vectors} $(\Hf^*,\Hg^*,\Hh^*)$ of tame level $N_{\Hf\Hg\Hh}=\lcm(N_\Hf,N_\Hg,N_\Hh)$ and then performs the construction of the $p$-adic $L$-function for this choice of test vectors, which grants some control on the nonvanishing of the local zeta-integrals at primes dividing $N_{\Hf\Hg\Hh}$ appearing in Ichino's formula (cf. \cite[theorem 1.1]{Ich2008}). In our applications finding the correct test vector will not be a problem, so the reader is invited to think of our generalized families $\Hg$ and $\Hh$ fixed above as test vectors for families of tame level dividing $M$.

\medskip
We show in section \ref{p-adicLFc} that the costruction of $\mathscr{L}^f_p(\Hf,\Hg,\Hh)$ can be extended to our more general setting.

\begin{prop}[cf. definition \ref{padicLfunction}, proposition \ref{padicLvalues} and proposition \ref{unbalinterpolation}]
Assume that the residual Galois representation $\overline{\V}_{\!\Hf}$ of the big Galois representation $\V_{\!\Hf}$ attached to $\Hf$ is absolutely irreducible and $p$-distinguished. Then there is an element $\mathscr{L}^f_p(\Hf,\Hg,\Hh)\in R_{\Hf\Hg\Hh}$ such that for every $\Hf$-unbalanced triple of meaningful weights $w=(x,y,z)$, the following formula holds:
\[
(\Ls_p^f(\Hf,\Hg,\Hh)(w))^2= \frac{L^*(\Hf_{\!x}\times\Hg_y\times\Hh_z,\frac{k+l+m-2}{2})}{\zeta_\Q(2)^2\cdot\Omega_{\Hf_{\!x}}^2}\cdot\mathscr{I}^{unb}_{w,\,p}\cdot\left(\prod_{\ell\mid M}\mathscr{I}^*_{w,\ell}\right)
\]
where:
\begin{itemize}
    \item [(i)] $L^*$ denotes the completed $L$-function (including the archimedean local factor);
    \item [(ii)] $\Omega_{\Hf_{\!x}}$ is a suitable period attached to $\Hf_{\!x}$, essentially given by its Petersson norm;        
    \item [(iii)] $\mathscr{I}^{unb}_{w,p}$ (resp. $\mathscr{I}^*_{w,\ell}$) is a suitable normalized local zeta integral at $p$ (resp. at $\ell$).
\end{itemize}
\end{prop}

\medskip

\begin{rmk}
\label{tripleproductliteraturecomp}
Here we group some observations elucidating the relations between our construction of $\Ls^f_p(\Hf,\Hg,\Hh)$ and the existing literature on the subject.
\begin{itemize}
    \item [(i)] As already pointed out, we adapt Hsieh's construction to our setting, following a method that essentially already appears in \cite[chapters 7 and 8]{Hi1993}. The theory of generalized $\Lambda$-adic forms developed in section \ref{generalized} allows us to simplify the construction. In particular, we show that the theory of ordinary parts carries over in this generalized setting (cf. proposition \ref{classicordR}) and thus we do not need to prove the equivalent of \cite[lemma 3.4]{Hsi2021}.
    \item [(ii)] The (only) novelty of our $p$-adic $L$-function consists in allowing $\Hg$ and $\Hh$ to be generalized families in the sense described above. It should be noted that this is not so surprising, since the so-called Panchishkin condition continues to hold \emph{in families} in the $\Hf$-unbalanced range when $\Hf$ is a Hida family and the families $\Hg$ and $\Hh$ are more general. We refer to \cite{Loe2023} for a more detailed explanation of why one should expect the existence of such $p$-adic $L$-functions. The $\Hf$-unbalanced triple product setting fits in this picture, as discussed briefly in \cite[section 5.3]{Loe2023}.
    \item [(iii)] In \cite{Fuk2022} the author provides a similar generalization of Hsieh's work to the case in which $\Hg$ and $\Hh$ are not necessarily Hida families. Yet, Fukunaga's notion of \emph{general $p$-adic families of modular forms} does not allow our generality for the rings of coefficients. Moreover, in the framework of \cite{Fuk2022} one cannot view the Fourier coefficients of such families as continuous/analytic functions on a suitable weight space in general.
    \item [(iv)] It would be interesting to find a way to extend our results to the case where $\Hf$ is a Coleman family (i.e., to the finite $p$-slope case), adapting the techniques of \cite{AI2021} (cf. also the recent preprint \cite{GPR2023}). We refer to \cite[section 5.5]{Loe2023} for a brief discussion on this point.
    \item [(v)] As already observed, we do not perform a general and careful \emph{level adjustment} as in \cite{Hsi2021}. It is clear that one could mimic Hsieh's recipe to achieve more generality in the construction.
    \item [(vi)] Thanks to the results of the preprint \cite{Mak2023} (which affords a canonical choice and an explicit interpolation formula for the so-called congruence number attached to a Hida family), the period $\Omega_{\Hf_{\!x}}$ becomes more explicit. We refer to remark \ref{canonicalcongrnumber} for more details.
\end{itemize}
\end{rmk}

\subsection{Factorization of triple product \texorpdfstring{$p$}{p}-adic \texorpdfstring{$L$}{L}-functions}
\label{factorizationsectionintro}
In the second part of the paper, we discuss some arithmetic applications in the setting the we now describe.

\medskip
Assume that $p\geq 5$ and let $\Hf$ be a Hida family of tame level $N_\Hf$ with trivial tame character. Fix $K/\Q$ a quadratic imaginary field of odd discriminant $-d_K$ and two ray class characters $\eta_1$ and $\eta_2$ of $K$, that we can view as valued in $L$.

\medskip
The following assumptions are in force:
\begin{itemize}
    \item [(A)] $p$ is inert in $K$;
    \item [(B)] $N_\Hf$ is squarefree, coprime to the discriminant of $K$ and with an even number of prime divisors which are inert in $K$ (\emph{Heegner hypothesis});
    \item [(C)] $\eta_i$ has conductor $c\,p^r\OO_K$, with $r\geq 1$ and $c\in\Z_{\geq 1}$, $(c\,, p\cdot d_K\!\cdot\! N_\Hf)=1$, $c$ not divisible by primes inert in $K$.
    \item [(D)] $\eta_1$ and $\eta_2$ are not induced by Dirichlet characters and the central characters of $\eta_1$ and $\eta_2$ are inverse to each other, so that $\varphi=\eta_1\eta_2$ and $\psi=\eta_1\eta_2^\sigma$ are ring class characters of $K$ (here $\langle \sigma\rangle =\Gal(K/\Q)$).
\end{itemize}

\medskip
A classical theorem of Hecke and Shimura attaches to the character $\eta_1$ (resp. $\eta_2$) a cuspidal newform $g$ (resp. $h$) of weight $1$, namely the theta series attached to $\eta_1$ (resp. $\eta_2$). In section \ref{CMfamilies} we describe how to realize $g$ (resp. $h$) as the weight $1$ specialization of a $p$-adic family $\Hg$ (resp. $\Hh$) of theta series of tame level $d_K$. Note that our notion of generalized $\Lambda$-adic form is taylored to include families such as $\Hg$ and $\Hh$ as non-trivial examples and that the specializations of $\Hg$ (resp. $\Hh$) will always be supercuspidal at $p$ (hence of infinite $p$-slope).

\medskip
After fixing a choice of test vectors $\Hg^*$ (resp. $\Hh^*$) of tame level $N_\Hf\!\cdot d_K\cdot c^2$, in section \ref{factorizationsection} we define an improved version $\mathcal{L}_p^f(\Hf,\Hg,\Hh)$ of $\Ls_p^f(\Hf,\Hg^*,\Hh^*)$, satisfying a simplified interpolation property. This relies on Hsieh's computations of local zeta integrals (and on Fukunaga's generalizations of Hsieh's results in \cite{Fuk2022}).

\medskip
Let $H_n$ denote the ring class field of $K$ of conductor $cp^n$ for every $n\in\Z_{\geq 0}$ and let $H_\infty$ be the union of all the $H_n$'s. Let $\mathscr{G}_\infty:=\Gal(H_\infty/K)$. We can identify the maximal $\Z_p$-free quotient $\Gamma^-$ of $\mathscr{G}_\infty$ with the Galois group of the anticyclotomic $\Z_p$-extension of $K$ and there is an exact sequence $0\to\Delta_c\to\mathscr{G}_\infty\to\Gamma^-\to 0$ of abelian groups with $\Delta_c$ a finite group and $\Gamma^-\cong\Z_p$. We fix a non-canonical isomorphism $\mathscr{G}_\infty\cong\Delta_c\times\Gamma^-$ once and for all.

Then $\varphi$ (resp. $\psi$) factors through $\mathscr{G}_\infty$ and we write it as $(\varphi_t,\varphi^-)$ (resp. $(\psi_t,\psi^-)$) according to the fixed isomorphism $\mathscr{G}_\infty\cong\Delta_c\times\Gamma^-$.

\medskip
For $k\in\Z_{\geq 2}\cap 2\Z$, let $\mathfrak{X}^{\mathrm{crit}}_{p,k}$ denote the set of continuous characters $\hat{\nu}:\Gamma^-\to\C_p^{\times}$ such that the associated algebraic Hecke character $\nu:\Ad_K^\times/K^\times\to\C^\times$ has infinity type $(j,-j)$ with $|j|<k/2$. 

\medskip
The main result of section \ref{factorizationsection} is the following factorization theorem for the \emph{anticyclotomic projection} $\mathcal{L}_{p,ac}^f(\Hf,\Hg,\Hh)$ of $\mathcal{L}_p^f(\Hf,\Hg,\Hh)$ (cf. definition \ref{anticycprojdefi}). This factorization is a counterpart of \cite[proposition 8.1]{Hsi2021} (which assumes $p$ split in $K$) and an upgrade of \cite[theorem 3.1]{BSVast2} to the case of Hecke characters with non-trivial $p$-part.

\begin{teo}[cf. theorem \ref{factthm}]
\label{factorizationthmintro}
In the above setting, it holds:
\[
\mathcal{L}^f_{p,ac}(\Hf,\Hg,\Hh)=\pm\mathscr{A}_{\Hf\Hg\Hh}\cdot\left(\varphi^-\left(\Theta_\infty^{\mathrm{Heeg}}(\Hf,\varphi_t)\right)~\hat{\otimes}~\psi^-\left(\Theta_\infty^{\mathrm{Heeg}}(\Hf,\psi_t)\right)\right)\,.
\]
\end{teo}

This equality takes place in the ring
\[
\mathcal{R}^-=\left(R_{\Gamma^-}\hat{\otimes}_\Lambda R_{\Gamma^-}\right)[1/p]\,,\quad\text{where}\quad R_{\Gamma^-}:=\Lambda\hat{\otimes}_{\OO_L}\OO_L[\![\Gamma^-]\!]
\]
and the notation is as follows.
\begin{itemize}
    \item [(i)] $\Theta_\infty^{\mathrm{Heeg}}(\Hf,\varphi_t)\in R_{\Gamma^-}$ (resp. $\Theta_\infty^{\mathrm{Heeg}}(\Hf,\psi_t)\in R_{\Gamma^-}$) is (a slight generalizations of) the so-called \emph{big theta element} constructed by Castella-Longo in \cite{CL2016}, building up on works by Bertolini-Darmon (cf. \cite{BD1996},\cite{BD1998},\cite{BD2007}) and Chida-Hsieh (cf. \cite{CH2018}). These $p$-adic $L$-functions interpolate (the square root of the algebraic part of) the special values $L(\Hf_{\!k}/K,\varphi_t\nu,k/2)$ (resp. $L(\Hf_{\!k}/K,\psi_t\nu,k/2)$) for $k\in\Z_{\geq 2}$ even and $\hat{\nu}\in\mathfrak{X}^{\mathrm{crit}}_{p,k}$.
    \item [(ii)] $\varphi^-(\tau)$ (resp. $\psi^-(\tau)$) for $\tau\in R_{\Gamma^-}$ denotes the image of the element $\tau$ via the $\OO_L$-linear automorphism of $R_{\Gamma^-}$ uniquely determined by the identity on $\Lambda$ and the assignment $[\gamma]\mapsto\varphi^-(\gamma)[\gamma]$ (resp. $[\gamma]\mapsto\psi^-(\gamma)[\gamma]$) on group-like elements on $\OO_L[\![\Gamma^-]\!]$.
    \item [(iii)] The element $\mathscr{A}_{\Hf\Hg\Hh}\in\mathcal{R}^-$ is defined in proposition \ref{adjustingfactor} and satisfies the crucial property that, for all $\hat{\nu},\hat{\mu}\in\mathfrak{X}^\mathrm{crit}_{p,2}$, $
\mathscr{A}_{\Hf\Hg\Hh}(2,\hat{\nu},\hat{\mu})\neq 0$.
\end{itemize}

\medskip
The proof of theorem \ref{factorizationthmintro} follows from the decomposition arising in our setting at the level of Galois representations (cf. lemma \ref{Galoisfact}) and from a careful comparison of the Euler factors at $p$ (or $p$-adic multipliers) appearing in the interpolation formulae for the various $p$-adic $L$-functions. In particular, this requires an explicit computation of the normalized local zeta integral at $p$ (denoted above by $\mathscr{I}^{unb}_{w,\,p}$), carried out in proposition \ref{localfactp}.

\subsection{\texorpdfstring{$p$}{p}-adic formulas for Heegner points} In section \ref{applicationssection} we apply theorem \ref{factorizationthmintro} to the study and the construction of Heegner points on elliptic curves. In what follows, we keep the notation as above and we let $E/\Q$ be an elliptic curve with multiplicative reduction at $p$. Let $f_E\in S_2(\Gamma_0(N_E))$ be the cuspidal newform of level $N_E$ attached to $E$ via modularity. Note that this implies that $N_E=p\cdot N_E^\circ$ with $p\nmid N_E^\circ$. Assume now that $\Hf$ denotes the unique primitive Hida family in $\SSS^\ord(N_E^\circ,\mathbbm{1},\Lambda)$ of tame level $N_E^\circ$ and trivial tame character, such that $\Hf_{\!2}=f_E$. 

We also impose an extra condition on the characters $\eta_1,\eta_2$ (cf. assumption \ref{quadraticphi}):
\begin{itemize}
    \item [(E)] $\varphi=\eta_1\eta_2$ has conductor prime to $p$ and $\psi=\eta_1\eta_2^\sigma$ has non-trivial anticyclotomic part (i.e., $\psi^-$ is non-trivial).
\end{itemize}

\medskip
In particular it follows that $\varphi^-$ is trivial and that we can identify $\varphi=\varphi_t$ as a character of the finite group $\Delta_c$. Let $H_\varphi$ denote the abelian extension of $K$ cut out by $\varphi$ and observe that $p$ splits completely in $H_\varphi$.

\medskip
Upon fixing a primitive Heegner point $P\in E(H_\varphi)\otimes\Q$ and setting $\alpha:=a_p(E)\in\{\pm 1\}$, one can define:
\begin{align*}
P_{\varphi}&:=\sum_{\sigma\in\Gal(H_{\varphi}/K)}\varphi(\sigma)^{-1} P^\sigma \in (E(H_{\varphi})\otimes\Q)^{\varphi}\\
P^{\pm}_{\varphi,\alpha}&:=P_\varphi\pm \alpha\cdot P_\varphi^{\Frob_\p}\in E(H_{\varphi})\otimes\Q\,.    
\end{align*}
One can show that $P^{\pm}_{\varphi,\alpha}$ does not depend on the choice of prime $\p$ of $H_\varphi$ above $p$. In what follows we fix the choice induced by our fixed embedding $\iota_p:\bar{\Q}\hookrightarrow\bar{\Q}_p$ and we view the points $P_\varphi$ and $P^{\pm}_{\varphi,\alpha}$ as elements of $E(\Q_{p^2})\otimes\Q$ under such an embedding.

\medskip
As $E$ has multiplicative reduction at $p$, we can take advantage of Tate's parametrization of $E$ to define a logarithm $\log_E: E(\Q_{p^2})\otimes\Q\to\Q_{p^2}$ at the level of $\Q_{p^2}$-rational points.

\medskip
Relying on theorem \ref{factorizationthmintro} and on previous results by Bertolini-Darmon (cf. \cite{BD1998} and \cite{BD2007}), we deduce the results summarized in the following statement.

\medskip
\begin{prop}[cf. corollaries \ref{derivkcor}, \ref{derivTcor} and \ref{infiniteorderiff}]
In the above setting, assume moreover that $L(E/K,\psi,1)\neq 0$. Then the restriction $\mathcal{L}_p^f(\Hf,g,h)$ of $\mathcal{L}_p^f(\Hf,\Hg,\Hh)$ to the \emph{line} $(k,1,1)$ vanishes at $k=2$ and
\[
\frac{d}{dk}\mathcal{L}_p^f(\Hf,g,h)_{|k=2}=\frac{c_E}{2}\cdot\log_E(P^+_{\varphi,\alpha})
\]
for some explicit constant $c_E\in\bar{\Q}_p^\times$. 

Similarly, the restriction $\mathcal{L}_{p,ac}^f(f_E,\Hg\Hh)$ of $\mathcal{L}_{p,ac}^f(\Hf,\Hg,\Hh)$ to the \emph{line}
$(2,\hat{\nu},\hat{\nu})$ vanishes at $\hat{\nu}=1$ (the trivial character) and
\[
\frac{d}{d\hat{\nu}}\mathcal{L}_{p,ac}^f(f_E,\Hg\Hh)_{|\hat{\nu}=1}= c_E\cdot\log_E(P^-_{\varphi,\alpha})
\]
for the same constant $c_E$. 

In particular, if $\varphi$ is a quadratic (or genus) character, the following are equivalent:
\begin{itemize}
    \item [(i)] 
\[
\left(\frac{d}{dk}\mathcal{L}_p^f(\Hf,g,h)_{|k=2}\,,\frac{d}{d\hat{\nu}}\mathcal{L}_{p,ac}^f(f_E,\Hg\Hh)_{|\hat{\nu}=1}\right)\neq(0,0)
\]
\item [(ii)] The point $P_\varphi$ is of infinite order.
\end{itemize}
\end{prop}

\begin{rmk}
In \cite{BSVast2} (cf. also \cite{DR2023}) the authors study a setting similar to ours, but require the characters $\eta_1$ and $\eta_2$ to have conductor coprime to $p$. As a consequence, the order of vanishing of the restriction $\mathcal{L}_p^f(\Hf,g,h)$ to the line $(k,1,1)$ of the corresponding triple product $p$-adic $L$-function is at least $2$. From a factorization in the style of theorem \ref{factorizationthmintro}, they deduce a formula for the second derivative of $\mathcal{L}_p^f(\Hf,\Hg,\Hh)$ at $k=2$ in terms of the product of logarithms of two Heegner points (respectively related to the characters that we denoted $\varphi$ and $\psi$). Our construction allows instead to pin down a single Heegner point from the study of $\mathcal{L}_p^f(\Hf,\Hg,\Hh)$ around the triple of weights $(2,1,1)$.    
\end{rmk}

\subsection*{Notation and conventions}
\medskip
If $F$ is any field, we denote by $G_F$ the absolute Galois group of $F$ (defined after fixing a suitable separable closure) and we denote $F^{ab}$ the maximal abelian extension of $F$ (inside such a separable closure).

\medskip
If $\Gamma$ is a profinite group and $R$ is a topological ring we denote by $R[\![\Gamma]\!]$ the completed group algebra with coefficients in $R$ (with the profinite topology) and we write $[\gamma]$ for $\gamma\in\Gamma$ to denote the corresponding group element in the ring $R[\![\Gamma]\!]$.

\medskip
We denote by $\Ad$ the ring of adéles of $\Q$ and if $B$ is a finite separable $\Q$-algebra we let $\Ad_B:=\Ad\otimes_\Q B$ denote the corresponding ring of adéles of $B$.

\medskip
For every number field $E$, we let the Artin reciprocity map
\[
\mathrm{rec}_E:\Ad_E^\times/E^\times\to\Gal(E^{ab}/E)
\]
to be \emph{arithmetically normalized}, i.e., if $v$ is a finite place of $E$ the compatible local Artin reciprocity map
\[
\mathrm{rec}_{E_v}:E_v^\times\to D_v\cong\Gal(E_v^{ab}/E_v)
\]
is the unique map such that for every uniformizer $\pi$ of $E_v$ it holds that $\mathrm{rec}_{E_v}(\pi)$ acts as the Frobenius morphism on the maximal unramified extension of $E_v$ (inside $E_v^{ab}$). We write $\Frob_v$ to denote an arithmetic Frobenius element at the place $v$ in $G_E$.

\medskip
If $K$ is a quadratic imaginary field and $\eta:G_K\to R^\times$ (here $R$ can be any ring) is a character, we let $\eta^\sigma$ to denote the conjugate of $\eta$, i.e., $\eta^\sigma(\gamma)=\eta(\sigma \gamma\sigma ^{-1})$ for $\gamma\in G_K$, where $\sigma\in G_K$ is any element such that $\sigma|_K$ generates $\Gal(K/\Q)$ (one possible explicit choice for $\sigma$ is the complex conjugation induced by the fixed embdedding $\iota_\infty$).

If $\chi:\Ad_K^\times/K^\times\to\C^\times$ is an algebraic Hecke character of $K$, we say that $\chi$ has $\infty$-type $(a,b)$ if for all $z\in\C^\times$ it holds $\chi(z\otimes 1)=z^{-a}\overline{z}^{-b}$.

\medskip
Given a smooth function $f$ on the upper-half plane $\HH:=\{ \tau\in\C\mid \Imm(\tau)>0\}$ and $\omega=\begin{psmallmatrix}a & b\\ c & d\end{psmallmatrix}\in\GL_2(\R)^+$ (invertible $2\times 2$ matrices with positive determinant) and $k\in\Z$, we set
\[
f|_k\omega(\tau):=\det(\omega)^{k/2}\cdot (c\tau +d)^{-k}\cdot f\left(\tfrac{a\tau+b}{c\tau+d}\right)\qquad \tau\in\HH
\]

\medskip
If $\Gamma\subseteq\SL_2(\Z)$ is a congruence subgroup and $k\in\Z_{\geq 1}$, we let $M_k(\Gamma)$ (resp. $S_k(\Gamma)$) be the $\C$-vector space of (holomorphic) modular forms (resp. cusp forms) of weight $k$ and level $\Gamma$. For $\Gamma=\Gamma_1(N)$ for some $N\geq 1$ and $\chi$ a Dirichlet character modulo $N$, we let $M_k(N,\chi)$ (resp. $S_k(N,\chi)$) denote the spaces of modular forms (resp. cusp forms) of weight $k$, level $\Gamma_1(N)$ and nebentypus $\chi$. Unless otherwise specified, we refer to \cite{Miy1989} for the all the basic facts concerning the analytic theory of modular forms which are mentioned freely without proof.


\subsection*{Acknowledgements}
This work contains parts of the author's PhD thesis written under the supervision of Massimo Bertolini. We would like to thank him for suggesting this project and for his precious advice. We thank the anonymous referee for their corrections, comments and suggestions. The author was funded by the DFG Graduiertenkolleg 2553 during his doctoral studies at Universität Duisburg-Essen and is currently funded by the Simons Collaboration on Perfection in Algebra, Geometry and Topology as postdoctoral researcher at CNRS. 

\section{Generalized \texorpdfstring{$\Lambda$}{lambda}-adic modular forms and ordinary projection}
\label{generalized}
In this section, we define a generalized notion of $\Lambda$-adic forms and we extend Hida's theory of the ordinary projector to this setting.
\subsection{First definition and examples}
Let $L$ be (as in the introduction) a finite extension of $\Q_p$, with ring of integers $\OO_L$, uniformizer $\varpi_L$ and residue field $\F_L:=\OO_L/\varpi_L\OO_L$.

Recall that $\Lambda:=\OO_L[\![1+p\Z_p]\!]$ is the completed group algebra for the profinite group $1+p\Z_p$. It is a complete local $\OO_L$ algebra of Krull dimension $2$, with maximal ideal $\mmm_\Lambda=(\varpi_L, T)$ and residue field $\F_L$. We fix once and for all the isomorphism
\[
\Lambda\cong\OO_L[\![T]\!]  
\]
uniquely determined by sending $[1+p]\mapsto 1+T$ and sometimes we write $\Lambda$ to denote directly $\OO_L[\![T]\!]$ via this identification.

\medskip
In this section, we will denote by $(R,\varphi)$ a complete local noetherian $\Lambda$-algebra (here we also mean that $\varphi:\Lambda\to R$ is a continuous local homomorphism of $\OO_L$-algebras) with maximal ideal $\mmm_R$ (also denoted $\mmm$ when it is clear from the context) and residue field $R/\mmm_R$ isomorphic to $\F_L$. We let $\CCC_{\Lambda}$ to be the category of such $\Lambda$-algebras, with arrows given by (continuous) homomorphisms of $\Lambda$-algebras. Similarly we have a category $\CCC_{\OO_L}$ and viewing $\Lambda$ as $\OO_L$-algebra in the obvious way, we get a functor $\CCC_{\Lambda}\to\CCC_{\OO_L}$ by pullback.

Sometimes we just write $R$ instead of $(R,\varphi)$ to simplify the notation, although the structure morphisms are going to play an important role in what follows.

\begin{defi}
For $R\in\CCC_{\Lambda}$ and any complete subring $\OO_L\subseteq A\subseteq\C_p$, we write
\[
\W_R(A):=\Hom^{cont}_{\OO_L-alg}(R,A)\,,
\]
endowed with the topology of uniform convergence on compact sets (which is essentially the $p$-adic topology). The elements of $\W_R(A)$ will be called ($A$-valued) $R$-weights (or $R$-specializations).
\end{defi}

\begin{rmk}
\label{intrmk}
Let $L'$ be a finite extension of $L$ inside $\C_p$ with ring of integers $\OO_{L'}$. Then for every $w\in\Hom^{cont}_{\OO_L-alg}(R,L')$ it holds $\OO_L\subseteq w(R)\subseteq L'$, but $w(R)$ cannot be a field. This forces $w(R)\subseteq \OO_{L'}$, so that we can identify $\W_R(L')=\W_R(\OO_{L'})=\Hom_{\CCC_{\OO_L}}(R,\OO_{L'})$ in our setting.
\end{rmk}

\medskip
We fix an embedding $\Z_p\hookrightarrow\W_\Lambda(L)$, given by sending $k\in\Z_p$ to the unique $\OO_L$-algebra homomorphism mapping $T\mapsto (1+p)^k-1$.

\begin{defi}
An element $w\in\W_\Lambda(\C_p)$ is an \emph{arithmetic} weight if it is uniquely determined by the assignment $T\mapsto\varepsilon(1+p)\cdot(1+p)^k-1$, where $k\in\Z_{\geq 1}$ and $\varepsilon:1+p\Z_p\to\mu_{p^\infty}(C_p^\times)$ is a finite order character. In this case we write $w=(k,\varepsilon)$ and we denote the set of arithmetic weights by $\W_\Lambda^{ar}$.

We say that $w=(k,\varepsilon)$ is \emph{classical} if $k\geq 2$ and we denote the set of classical weights by $\W_\Lambda^{cl}$. Clearly $\Z_p\cap\W_\Lambda^{cl}=\Z_{\geq 2}\subset\W_\Lambda^{cl}$ via the embedding $\Z_p\hookrightarrow\W_\Lambda(L)$.
\end{defi}

\begin{defi}
Let $(R,\varphi)\in\CCC_\Lambda$. We define the set of classical $R$-weights as
\[
\W_R^{cl}:=\{w\in\W_R(\C_p)\mid w\circ\varphi\in\W_\Lambda^{cl}\}
\]
and the set of integral classical $R$-weights as
\[
\W^{cl}_{R,\Z}:=\{w\in\W_R(\C_p)\mid w\circ\varphi\in\Z_{\geq 2}\}.
\]
For every $w\in\W_R^{cl}$ we define $(k_w,\varepsilon_w):=w\circ\varphi$ and, if $w\circ\varphi\in\W^{cl}_{R,\Z}$, we simply write $w\circ\varphi=k_w$. For any subset $V\subset\W_R(\C_p)$ we set $\varphi^*(V)=\{w\circ\varphi\mid w\in V\}$.  
\end{defi}

\begin{defi}
\label{decentsetwt}
We say that a subset $\Omega\subseteq\W^{cl}_{R,\Z}$ is \emph{$(\Lambda,R)$-admissible} if the following conditions are satisfied:
\begin{itemize}
    \item [(i)] the closure of $\varphi^*(\Omega)$ inside $\Z_p\subseteq\W_\Lambda(L)$ contains a non-empty open subset of $\Z_p$;
    \item [(ii)] $\bigcap_{w\in\Omega}\Ker(w)=\{0\}$.
\end{itemize}
\end{defi}

\medskip
We will need the following result later.
\begin{lemma}
\label{compliso}
Let $R\in\CCC_\Lambda$ and let $\CS\subseteq\W^{cl}_R$ be a countable infinite set. Let $\BB$ denote the set of ideals in $R$ that can be written as a finite intersection of pairwise different primes of $R$ of the form $\q=\Ker(w)$ for $w\in\CS$. For every $J\in\BB$, consider $R/J$ with the quotient topology. Let $I=\bigcap_{w\in\CS} \Ker(w)$ and consider $R/I$ with the quotient topology. Then the natural map $R/I\to\varprojlim_{J\in\BB} R/J$ induces an isomorphism of topological rings $R/I\cong\varprojlim_{J\in\BB} R/J$.
\begin{proof}
For every $J\in\BB$, $R/J$ is a complete noetherian local ring with maximal ideal $\mmm_R/I$. Note that the quotient topology and the $\mmm_R/J$-adic topology on $R/J$ coincide and that the natural projection $R\to R/J$ is open and continuous (the same applies to $R/I$).

We claim that such topology on $R/J$ is the same as the $\varpi_L$-adic topology. It is clear that for every $n\geq 1$ it holds that $(\varpi_L^nR+J)/J\subseteq(\mmm_R^n+J)/J$. We are left to show that, for every $n\geq 1$, $R/(J,\varpi_L^n)$ is a quotient of $R/\mmm_R^m$ for $m\gg 1$ (in particular it is a finite ring). Indeed writing $J=\q_1\cap\dots\cap\q_s$ one checks that
\[
\sqrt{(J,\varpi_L^n)}=\sqrt{\bigcap_{i=1}^s(\q_i,\varpi_L^n)}=\bigcap_{i=1}^s\sqrt{(\q_i,\varpi_L^n)}=\bigcap_{i=1}^s\sqrt{(\q_i,\varpi_L)}=\mmm_R.
\]
The first equality follows from $\left(\bigcap_{i=1}^s(\q_i,\varpi_L^n)\right)^s\subseteq (J,\varpi_L^n)\subseteq\bigcap_{i=1}^s(\q_i,\varpi_L^n)$.
The second and the third equalities are obvious. The last one follows from the fact that $\sqrt{(\q_i,\varpi_L)}=\mmm_R$ for all $i=1,\dots,s$, since $R/\q_i$ is (algebraically isomorphic to) a finite extension of $\OO_L$ inside $\bar{Q}_p$ and $R/\mmm_R=\F_L$ by assumption. In particular it follows that $\mmm_R^m\subseteq (J,\varpi_L^n)$ for some $m\geq 1$ large enough, proving our claim. Hence we have natural topological isomorphisms for all $J\in\BB$
\[
R/J\cong\varprojlim_n R/(J,\varpi_L^n).
\]

Arguing as above it also follows that a fundamental system of open neighbourhoods of $0$ in $R/I$ is given by the open ideals $\{(\varpi_L^n+J)/I)\}_{n\geq 1, J\in\BB}$.

\medskip
This shows that we can realize the natural map $R/I\to\varprojlim_{J\in\BB}R/J$ as a chain of topological isomorphisms
\[
R/I\cong\varprojlim_{J\in\BB,\,n}R/(J,\varpi_L^n)\cong\varprojlim_{J\in\BB}R/J\,,
\]
concluding the proof.
\end{proof}
\end{lemma}

\medskip
We are ready to give the key definition of this section.

\begin{defi}
Let $N\in\Z_{\geq 1}$ be an integer with $p\nmid N$, let $\chi$ be a Dirichlet character modulo $Np^t$ for some $t\in\Z_{\geq 1}$ with values in $\OO_L^{\times}$. We say that a \emph{generalized $\Lambda$-adic form} of tame level $N$ and \emph{branch character} $\chi$ is a couple $((R,\varphi),\Hxi)$ where:
\begin{itemize}
    \item [(i)] $(R,\varphi)$ is an object of $\CCC_\Lambda$, which is also flat as $\Lambda$-algebra and an integral domain,
    \item [(ii)] $\Hxi\in R[\![q]\!]$ is a formal $q$-expansion,
 \end{itemize}
such that the set of integral weights
\[
\Omega_{\Hxi,\Z}:=\{w\in\W^{cl}_{R,\Z}\mid \Hxi_{\!w}\in M_{k_w}(Np^t,\chi\omega^{2-k_w},\C_p)\}
\]
is $(\Lambda,R)$-admissible in the sense of definition \ref{decentsetwt}, where $\Hxi_{\!w}$ denotes the $q$-expansion obtained applying $w$ to the coefficients of $\Hxi$. We say that $((R,\varphi),\Hxi)$ is \emph{cuspidal} if, moreover, $\Hxi_{\!w}$ is cuspidal for all $w\in\Omega_{\Hxi,\Z}$.

Given a generalized $\Lambda$-adic form $((R,\varphi),\Hxi)$ and a $(\Lambda,R)$-admissible set of integral classical weights $\Omega\subseteq\Omega_{\Hxi,\Z}$, we say that $((R,\varphi),\Hxi)$ is \emph{$\Omega$-compatible}. Often we shorten the notation and we simply write $\Hxi$ to denote the $\Lambda$-adic form $((R,\varphi),\Hxi)$.
\end{defi}

\begin{defi}
Given a generalized $\Lambda$-adic form of tame level $N$ and character $\chi$ with coefficients in $(R,\varphi)$, we set
\[
\Omega_{\Hxi}:=\{w\in\W^{ar}_R\mid\Hxi_{\!w}\in M_{k_w}(Np^{e_w},\chi\omega^{2-k_w}\varepsilon_w,\C_p)\}
\]
where the exponent $e_w\geq 1$ depends on the $p$-part of $\chi$ and on $w$.
\end{defi}

\begin{defi}
We let $\M_\Omega(N,\chi,(R,\varphi))$ (respectively $\SSS_\Omega(N,\chi,(R,\varphi))$) denote the $R$-modules of generalized $\Lambda$-adic forms (resp. cuspidal generalized $\Lambda$-adic forms) of level $N$ and character $\chi$, with coefficients in $(R,\varphi)$ and $\Omega$-compatible (where $\Omega$ is a $(\Lambda,R)$-admissible set of classical integral $R$-weights). When all the inputs are clear from the context (or when it is not necessary to specify them) we simply write $\M$ and $\SSS$ to denote such $R$-modules, which we view as submodules of $R[\![q]\!]$ in the obvious way. We endow all such $R$-modules with the $\mmm$-adic topology. 
\end{defi}

\begin{rmk}
\label{R[[q]]}
The noetherianity of $R$ implies that $R[\![q]\!]$ is $\mmm$-adically separated and complete.
\end{rmk}

\begin{rmk}
\label{lambdaadicHeckeop}
On $\M=\M_\Omega(N,\chi,(R,\varphi))$ and $\SSS=\SSS_\Omega(N,\chi,(R,\varphi))$) there is an action of Hecke operators $T_\ell$ for $\ell\nmid Np$ prime, $U_\ell$ for $\ell\mid N$ prime and $U_p$. Those operators can be defined directly on the $q$-expansions in such a way that the specialization maps are Hecke-equivariant morphisms. More precisely, there is a character $\langle\,\cdot\,\rangle_\Lambda:\Z_p^\times\to\Lambda^\times$ given by $\langle s\rangle = [s\cdot\omega^{-1}(s)]$. For $(R,\varphi)\in\CCC_\Lambda$ we then let $\langle\,\cdot\,\rangle_R:\Z_p^\times\to R^\times$ to be the composition of $\langle\,\cdot\,\rangle_\Lambda$ with $\varphi$. Then, for every $\Hxi=\sum_{n=0}a_n(\Hxi)q^n\in\M$ and for every prime $\ell\neq p$ the Hecke operator $T_\ell$ acts as follows
\[
T_\ell(\Hxi)=\sum_{n=0}^{+\infty}a_n(T_\ell(\Hxi))q^n\,,\quad\text{where}\quad a_n(T_\ell(\Hxi))=\sum_{d\mid (n,\ell)}\langle d\rangle_R\cdot \chi(d)d^{-1}a_{n\ell/d^2}(\Hxi),
\]
with the convention that $\chi(\ell)=0$ if $\ell\mid N$. If $\ell\mid N$, we write $U_\ell$ to denote the $T_\ell$ operator

\medskip
We are particularly interested about the $U_p$ operator, whose action on $q$-expansions is the familiar one:
\[
U_p\left(\sum_{n=0}^{+\infty}a_nq^n\right)=\sum_{n=0}^{+\infty}a_{np} q^n \,. 
\]

\medskip
We end this remark recalling the action of the $V_p$ operator on $q$-expansions, given by
\[
V_p\left(\sum_{n=0}^{+\infty}a_nq^n\right)=\sum_{n=0}^{+\infty}a_p q^{np} \,. 
\]
This operator will appear later in the paper. Recall that $U_p\circ V_p$ is the identity on $q$-expansions, while $1-V_p\circ U_p$ defines the so-called $p$-depletion operator.
\end{rmk}

\begin{defi}
Let $N$, $\chi$, $(R,\varphi)$ and $\Omega$ be as above. The notation $\T_\Omega(N,\chi,(R,\varphi))$ will denote the $R$-subalgebra of $\End_R(\SSS_\Omega(N,\chi,(R,\varphi))$ generated by the Hecke operators $T_\ell$ for $\ell\nmid Np$ prime, $U_\ell$ for $\ell\mid N$ prime and $U_p$. When all the inputs are clear from the context we simply write $\T$ or $\T_\Omega$ to denote such Hecke algebra.
\end{defi}

\begin{defi}
An element $\Hxi\in\M$ is called a \textbf{generalized $\Lambda$-adic eigenform} (of given tame level $N$, character, branch, coefficients) if it is a simultaneous eigenvector for the Hecke operators $T_\ell$ ($\ell\nmid Np$ prime) and for the Hecke operator $U_p$.
\end{defi}

\begin{ex}
\label{prodfamex}
Let $\Hxi_1\in\M_{\Omega_1}(N,\chi_1,R_1)$ and $\Hxi_2\in\M_{\Omega_2}(N,\chi_2, R_2)$. Set $R:=R_1\hat{\otimes}_{\OO_L} R_2$. If $\mmm_i\subset R_i$ denotes the respective maximal ideal for $i=1,2$, then recall that by definition
\[
R=\varprojlim_{m,n}\left(\tfrac{R_1}{\mmm_1^n}\otimes_{\OO_L}\tfrac{R_2}{\mmm_2^m}\right).
\]
$R$ is then identified with the $\tilde{\mmm}$-adic completion of $R_1\otimes_{\OO_L}R_2$ where
\[
\tilde{\mmm}=\mmm_1\otimes_{\OO_L} R_2 + R_1\otimes_{\OO_L}\mmm_2\subset R_1\otimes_{\OO_L}R_2
\]
is a maximal ideal of $R_1\otimes_{\OO_L}R_2$ such that $(R_1\otimes_{\OO_L}R_2)/\tilde{\mmm}\cong\F_L$ (thanks to our strict conditions on the residue fields of $R_1$ and $R_2$). 

\medskip
For every $a\in R_1, b\in R_2$ we let $a\hat{\otimes}b$ denote the image of $a\!\otimes\!b\in R_1\otimes_{\OO_L}R_1$ inside $R$ via the natural map.
We endow $R$ with the following canonical $\Lambda$-algebra structure $\varphi :\Lambda\to R$ uniquely determined by $\OO_L$-linearity and the assignment
\[
\varphi(T):=\varphi_1(T)\hat{\otimes}1+1\hat{\otimes}\varphi_2(T)+\varphi_1(T)\hat{\otimes}\varphi_2(T)
\]
where $\varphi_i$ are the structure morphisms for $R_i$, $i=1,2$ (notice that this is well-defined).

\medskip
We refer to \cite[section 0.7.7]{EGAI} for the needed properties of completed tensor products. In particular it follows that $R\in\CCC_\Lambda$ and $R$ is an integral domain. Note that $R$ is a flat $\Lambda$-algebra via $\varphi$. This can be seen easily factoring $\varphi$ as composition of flat morphisms as
\[
\Lambda\to\Lambda\hat{\otimes}_{\OO_L}\Lambda\xrightarrow{(\varphi_1,\varphi_2)}R,
\]
where the first arrow sends $T\mapsto T\hat{\otimes}1+1\hat{\otimes}T+T\hat{\otimes}T$.

\medskip
By the universal property of completed tensor product it follows that, for every complete subring $A$ of $\C_p$ containing $\OO_L$ ,$\W_R(A)=\W_{R_1}(A)\times\W_{R_2}(A)$ (also as topological spaces) and, by our definition of $\varphi$, it also follows that under this identification we get an inclusion
\[
\W^{cl}_{R_1,\Z}\times\W^{cl}_{R_2,\Z}\subset\W^{cl}_{R,\Z}
\]
such that
\[
k_{(w_1,w_2)}=(w_1,w_2)\circ\varphi=(w_1\circ\varphi_1)+(w_2\circ\varphi_2)=k_{w_1}+k_{w_2}
\]
Let $\Omega=\Omega_1\times\Omega_2$, viewed as a subset of $\W^{cl}_{R,\Z}$ as above. It is easy to see that $\Omega$ is $(\Lambda,R)$-admissible.
 
\medskip
It follows that $\Hxi_1\times\Hxi_2\in\M_\Omega(N,\chi_1\chi_2\omega^2,R)$, where as usual if
\[
\Hxi_1=\sum_{n=0}^{+\infty}a_nq^n,\qquad \Hxi_1=\sum_{n=0}^{+\infty}b_nq^n
\]
we let
\[
\Hxi_1\times\Hxi_2=\sum_{n=0}^{+\infty}\left(\sum_{j=0}^n a_j\hat{\otimes}b_{n-j}\right)q^n\in R[\![q]\!].
\]

Indeed it is clear that, for all $(w_1,w_2)\in\Omega_1\times\Omega_2$, it holds
\[
(\Hxi_1\times\Hxi_2)_{(w_1,w_2)}=\Hxi_{1,w_1}\times \Hxi_{2,w_2}\in M_{k_{(w_1,w_2)}}(Np^t,\chi_1\chi_2\cdot \omega^{4-k_{(w_1,w_2)}},\OO_L)
\]

\end{ex}

\subsection{The ordinary projector}
We want to check that also in our generalized setting one can attach to the operator $U_p$ an idempotent operator $e^\ord$ obtained as
\[
e^\ord=\lim_{n\to +\infty}U_p^{n!}
\]
where the limit is taken in the $\mmm$-adic topology.

\begin{lemma}
The $R$-modules $\M$ and $\SSS$ are $\mmm$-adically complete and separated for the $\mmm$-adic topology.
\begin{proof}
We only give the proof for $\M=\M_\Omega$ (the proof for $\SSS$ is identical). It is clear that $\M$ is $\mmm$-adically separated, being a submodule of $R[\![q]\!]$ (which is $\mmm$-adically complete and separated by remark \ref{R[[q]]}). An element $(\overline{\Hxi}_n)_{n\geq 1}\in\varprojlim_{n}\M/\mmm^n\M$ defines (by left exactness of $\varprojlim_n$) a unique element
\[
\Hxi\in R[\![q]\!]=\varprojlim_n R[\![q]\!]/\mmm^n R[\![q]\!]\,.
\]
If for every $n\geq 1$ we fix a lift $\Hxi_n\in\M$ of $\overline{\Hxi}_n$ we know that for every $w\in\Omega$ it holds $\Hxi_{n,w}\in M_{k_w}(Np^t,\chi\omega^{2-{k_w}},\OO_L)$ and by the continuity of the specializations and the fact that $M_{k_w}(Np^t,\chi\omega^{2-{k_w}},\OO_L)$ is a finite and free $\OO_L$-module (thus complete), we deduce that
\[
\Hxi_{\!w}=\lim_{n\to +\infty}\Hxi_{n,w}\in M_{k_w}(Np^t,\chi\omega^{2-{k_w}},\OO_L),
\]
so that indeed $\Hxi\in\M$ and the lemma follows.
\end{proof}
\end{lemma}

\begin{prop}
\label{propordproj}
There exists a unique ordinary projector $e^\ord\in\End_R(\M)$ attached to the Hecke operator $U_p$, such that
\begin{enumerate}
    \item [(i)] $e^\ord(\Hxi)=\lim_{n\to +\infty}U_p^{n!}(\Hxi)$ (limit taken in the $\mmm$-adic topology)
    \item [(ii)] $e^\ord$ and $U_p$ commute and the module $\M$ carries a $U_p$-stable decomposition $\M = e^\ord\M\oplus (1-e^\ord)\M$ where $U_p$ is bijective on $e^\ord\M$ and topologically nilpotent on $(1-e^\ord)\M$.
    \item [(iii)] $e^\ord$ commutes with $T_\ell$ for all $\ell\nmid Np$ and is compatible with every meaningful arithmetic specialization.
    \item [(iv)] the formation of $e^\ord$ is compatible with inclusions $\M_\Omega\subseteq\M_{\Omega'}$ induced by inclusions $\Omega'\subseteq\Omega$ of $(\Lambda,R)$-admissible sets of classical integral weights.
\end{enumerate}
The analogue assertions for $\SSS$ hold.
\begin{proof}
The proof of the proposition follows well-known ideas, so we just sketch the argument. We only give the proof for $\M=\M_\Omega$ (the proof for $\SSS$ is identical).

Condition (ii) of definition \ref{decentsetwt} implies that the specialization maps give rise to an inclusion
\begin{equation}
\label{lambdaadicinclusion}
\M_\Omega\hookrightarrow\prod_{w\in\Omega}M_{k_w}(Np^t,\chi\omega^{2-k_w},\OO_L)\,.
\end{equation}
Note that, since for each $w\in\Omega$ we have $\Ker(w)+\varpi_L R=\mmm$, the above inclusion is an embedding for the $\mmm$-adic topology on the LHS and the $\varpi_L$-adic topology on the RHS. It is well-known (cf. \cite[section 7.2]{Hi1993}) that on the RHS there exists an ordinary projector associated with the operator $U_p$ (where the limit $\lim_{n\to +\infty}U_p^{n!}$ is taken in the $\varpi_L$-adic topology). Since the action of $U_p$ commutes with specializations and $\M_\Omega$ is $\mmm$-adically separated and complete,  by restriction along \eqref{lambdaadicinclusion} we obtain an idempotent $e^\ord$ on $\M_\Omega$ which visibly satisfies all the required properties.
\end{proof}
\end{prop}

\medskip
We are then led to the following definition:

\begin{defi}
\label{genHidafam}
We say that a generalized eigenform $\Hxi\in\M_\Omega(N,\chi,R)$ (respectively $\Hxi\in\SSS_\Omega(N,\chi,R)$) is a \emph{generalized Hida family} (resp. a \emph{cuspidal generalized Hida family}) if $e^\ord(\Hxi)=\Hxi$.

We define the $R$-modules $\M^\ord_\Omega(N,\chi,R):=e^\ord(\M_\Omega(N,\chi,R))$ (resp. in the cuspidal case $\SSS^\ord_\Omega(N,\chi,R):=e^\ord(\SSS_\Omega(N,\chi,R))$) to be the submodules of $\M_\Omega(N,\chi,R)$ (resp. $\SSS_\Omega(N,\chi,R)$) of ordinary generalized $\Lambda$-adic forms. When the inputs are clear from the context we simply write $\M^\ord$ or $\M^\ord_\Omega$ (resp. $\SSS^\ord$ or $\SSS^\ord_\Omega$).

We let $\T^\ord_\Omega(N,\chi,R)$ to denote the $R$-subalgebra of $\End_R(\SSS^\ord_\Omega(N,\chi,R))$ generated by the Hecke operators $T_\ell$ for $\ell\nmid Np$ prime, $U_\ell$ for $\ell\mid N$ prime and $U_p$. When all the inputs are clear from the context, we simply write $\T^\ord$ or $\T^\ord_\Omega$ to denote such Hecke algebra.
\end{defi}

\begin{rmk}
Equivalently one could define generalized Hida families asking that every meaningful classical specialization is a $p$-ordinary eigenform in the usual sense.
\end{rmk}

\medskip
The following proposition shows that generalized Hida families are actually essentially the same as classical Hida families.

\begin{prop}
\label{classicordR}
For any $R\in\CCC_\Lambda$ which is $\Lambda$-flat and an integral domain and any $(\Lambda,R)$-admissible set of classical integral weights $\Omega$, the $R$-modules $\M^\ord_\Omega(N,\chi,R)$ (resp. $\SSS^\ord_\Omega(N,\chi,R)$) are free $R$-modules of finite rank. Moreover (assuming that $\chi$ takes values in $\OO_L^\times$), there are canonical isomorphisms
\[
\M^\ord(N,\chi,\Lambda)\otimes_\Lambda R\xrightarrow{\cong}\M^\ord_\Omega(N,\chi,R),\quad \SSS^\ord(N,\chi,\Lambda)\otimes_\Lambda R\xrightarrow{\cong}\SSS^\ord_\Omega(N,\chi,R)\,.
\]
\begin{proof}
We will omit the proof of the cuspidal case because the proof does not change. In this proof, we write $\M^\ord_\Lambda=\M^\ord(N,\chi,\Lambda)$ and $\M^\ord_R=\M^\ord(N,\chi,R)$ to simplify the notation.
In order to prove that $\M^\ord_R$ is $R$-free of finite rank we adapt Wiles's proof for classical Hida theory (cf. \cite[section 7.3]{Hi1993}). We recall the main ideas for the convenience of the reader. Let $M$ be a finite free $R$-submodule of $\M^\ord_R$, with $R$-basis $\{\Hxi_1,\dots,\Hxi_r\}$. Write
\[
\Hxi_i=\sum_{n=0}^{+\infty} a_n(\Hxi_i)q^n
\]
for $i=1,\dots,r$. Then there is a sequence of integers $0\leq n_1<n_2<\dots<n_r$ such that the $r\times r$ matrix $\left(a_{n_j}(\Hxi_i)\right)_{i,j,=1,\dots r}$ has non-zero determinant $d\in R$. Since by assumption $\cap_{w\in\Omega}\Ker(w)=(0)$, we deduce that there exists $w\in\Omega$ such that $d\neq 0\mod \Ker(w)$, so that the specializations $\{\Hxi_{1,w},\dots,\Hxi_{r,w}\}$ would still be $\OO_L[w]$-linearly independent in $M^\ord_{k_w}(Np^t,\chi\omega^{2-k_w},\OO_L[w])$. It is well-known (and established by Hida) that the rank of $M^\ord_{k_w}(Np^t,\chi\omega^{2-k_w},\OO_L[w])$ is independent on $w$ if $k_w\geq 3$. Hence there exists $r^*\in\Z_{\geq 0}$ such that $\M^\ord_R$ admits finite free $R$-submodules of rank $r^*$, but not of rank $r^*\!+1$. Assume now that $M$ is such a finite free $R$-submodule of $\M^\ord_R$ of rank $r^*$. One checks easily that, with the notation as above, $d\cdot\M^\ord_R\subseteq M$. Hence, by the noetherianity of $R$, it follows that $\M^\ord_R$ is finitely generated as $R$-module. In particular it is a compact $R$-module (equivalently a profinite $R$-module). The topological Nakayama's lemma (cf. \cite[lemma 3.2.6]{Hid2012} for instance) implies that $\M_R^\ord$ is generated by $r:=\dim_{\F_L}(\M_R^\ord/\mmm_R \M_R^\ord)$ elements (a lift of an $\F_L$-basis of $\M_R^\ord/\mmm_R \M_R^\ord$).

Now note that (using the flatness of $R$ over $\Lambda$) $\M^\ord_\Lambda\otimes_\Lambda R$ can be naturally seen as an $R$-free submodule of $\M^\ord_R$ of $R$-rank $r$. We define the quotient
\[
Q:=\frac{\M^\ord_R}{\M^\ord_\Lambda\otimes_\Lambda R}
\]
and we claim that $Q=0$. This would conclude the proof of the proposition, since it is well-known that $\M^\ord_\Lambda$ is a free $\Lambda$-module of rank $r^*$.

Picking $w\in\Omega$ with $k_w\geq 3$, one has $Q\otimes_R R/\Ker(w)=0$, since both $\M^\ord_\Lambda\otimes_\Lambda R$ and $\M^\ord_R$ project onto $M^\ord_{k_w}(Np^t,\chi\omega^{2-k_w},\OO_L[w])$ via $w$ (to see this one uses the trick of twisting with a suitable family of Eisenstein series, cf. \cite[pag. 199]{Hi1993}). Hence \textit{a fortiori} $Q\otimes_R R/\mmm= 0$ and, since also $Q$ is a profinite $R$-module, it follows again from the topological Nakayama's lemma that $Q=0$.
\end{proof}
\end{prop}

\begin{rmk}
\label{omegaindip}
Proposition \ref{classicordR} shows that the $R$-modules $\M^\ord_\Omega(N,\chi,R)$ (respectively  $\SSS^\ord_\Omega(N,\chi,R)$) actually does not depend on $\Omega$, so that in the ordinary setting we will omit the $(\Lambda,R)$-admissible set of weights from the notation from now on.
\end{rmk}
\section{The unbalanced triple product \texorpdfstring{$p$}{p}-adic \texorpdfstring{$L$}{L}-function}
\label{p-adicLFc}
In this section we carry out the construction of a generalized unbalanced triple product \texorpdfstring{$p$}{p}-adic \texorpdfstring{$L$}{L}-function, closely following the method appearing in \cite{Hsi2021}. Having defined the ordinary projector $e^\ord$ in wider generality and having proved proposition \ref{classicordR}, the construction simplifies slightly. For instance, we do not need the equivalent of \cite[lemma 3.4]{Hsi2021}.

\subsection{Remarks on the Atkin-Lehner involution}
Recall that given $\xi\in S_k(M,\chi)$, one has an Atkin-Lehner involution $w_M: S_k(M,\chi)\to S_k(M,\chi^{-1})$ given by $w_M(\xi)=\xi|_k\begin{psmallmatrix}
    0 & -1\\
    M & 0
\end{psmallmatrix}$. For our constructions we will need a $\Lambda$-adic version of the Atkin-Lehner involution. This entails considering more general Atkin-Lehner operators.

\medskip
Let $N$ be a positive integer coprime to $p$ and $t\in\Z_{\geq 1}$. If $d$ is an integer coprime to $Np$, we write $\langle d\rangle = \langle a;b\rangle$ for the diamond operator corresponding to $d\in(\Z/Np^{t}\Z)^\times$, where the convention is that $d\equiv a\mod N$ and $d\equiv b\mod p^t$.

For $\xi\in S_k(Np^t,\chi)$ we define the Atkin-Lehner operator $w_N$ on $\xi$ as 
\begin{equation}
\label{wnoperator}
    w_N(\xi):=\langle 1 ; N\rangle (\xi|_k\,\omega_N)\qquad \omega_N:=\omega_{N,\,p^t}:=\begin{pmatrix}N & -1 \\ Np^tc & Nd\end{pmatrix}\,,
\end{equation}
where we require that $\det(\omega_N)=N$. Write $\chi=\chi_{p^t}\chi_{N}$ in a unique way for $\chi_{p^t}$ a character modulo $p^t$ and $\chi_N$ a character modulo $N$. 

Then (cf. \cite[\S 1]{AL1978}, where they define an operator which is the \emph{inverse} of ours) $w_N$ is an operator
\[
w_N:S_k(Np^t,\chi)\to S_k(Np^t,\overline{\chi_N}\chi_{p^t})
\]
such that for all primes $\ell\nmid N$ it holds that $w_N\circ T_\ell=\chi_N(\ell)(T_\ell\circ w_N)$ and (when $t\geq 1$) that $w_N\circ U_p=\chi_N(p)(U_p\circ w_N)$. One can also check that if $s>r\geq 0$, the action of $w_N$ on $S_k(\Gamma_1(Np^r))$ is the restriction of the action of $w_N$ on $S_k(\Gamma_1(Np^s))$, by our choice of the matrices $\omega_{N,\,p^t}$, so that it makes sense to drop $p^t$ from the notation.

\medskip
In particular, if $\xi\in S_k(Np^t,\chi)$ is a normalized newform, then $w_N(\xi)=\lambda_N(\xi)\cdot \breve{\xi}$ where $\lambda_N(\xi)$ is an algebraic number of complex absolute value $1$ (a so called \emph{pseudo-eigenvalue}) and $\breve{\xi}$ is a normalized newform such that if
\[
\xi=\sum_{n=1}^{+\infty}a_n q^n\qquad \breve{\xi}=\sum_{n=1}^{+\infty}b_n q^n
\]
then
\[
b_\ell=\begin{cases}
    \overline{\chi_N}(\ell)a_\ell & \text{ if }\ell\nmid N \\
    \chi_{p^t}(\ell)a_\ell & \text{ if }\ell\mid N.
    \end{cases}
\]

Moreover if $\xi\in S_k(N,\chi)$ is a $p$-ordinary newform with $k\geq 2$ and $\xi_\alpha\in S_k(Np,\chi)$ is its ordinary $p$-stabilisation, then $\lambda_N(\xi)^{-1}\cdot w_N(\xi_\alpha)$ coincides with the ordinary $p$-stabilisation of the newform $\breve{\xi}$, so we will write
\[
\breve{\xi_\alpha}:=\lambda_N(\xi)^{-1}\cdot w_N(\xi_\alpha).
\]
Note that in this case it is well-known that $\breve{\xi}$ is the modular form obtained applying complex conjugation to the Fourier coefficients of $\xi$.

\bigskip
Now let
\[
\Hxi=\sum_{n=1}^{+\infty}a_n(\Hxi)q^n\in \SSS^\ord(N_{\Hxi},\chi_{\Hxi},\Lambda_{\Hxi})
\]
be a classical new Hida family of tame level $N_{\Hxi}$ with character $\chi_{\Hxi}$ of conductor dividing $N_{\Hxi}\cdot p$, i.e., the classical specializations at integral weights of $\Hxi$ are either newforms of level $N_{\Hxi}\cdot p$ or ordinary $p$-stabilizations of newforms of level $N_{\Hxi}$. Here $\Lambda_{\Hxi}$ is a finite flat $\Lambda$-algebra in $\CCC_\Lambda$ and we assume that $L$ contains a primitive $N_{\Hxi}$-th root of unity. We require that $\Hxi$ is normalized (i.e., $a_1(\Hxi)=1$). Note that we can omit the admissible set of integral classical weights in the notation here, since classical Hida theory shows that for classical Hida families it always happens $\Omega_{\Hxi,\Z}=\W^{cl}_{\Lambda_{\Hxi},\Z}$.

\medskip
Following \cite[section 3.3]{Hsi2021}, there is a unique new Hida family $\breve{\Hxi\,}\in\SSS^\ord(N_{\Hxi} ,\chi^{-1}_{\Hxi},\Lambda_{\Hxi})$ which is characterised by the fact that, for all $x\in \W_{\Lambda_{\Hxi}}^{cl}$
\[
 (\breve{\Hxi\,})_x=\breve{(\Hxi_{x})}=\lambda_N(\Hxi_{x})^{-1}\cdot w_N(\Hxi_{x}).
\]

\medskip
\subsection[\texorpdfstring{Construction of the $p$-adic $L$-function}{Construction of the p-adic L-function}]{Construction of the $p$-adic $L$-function}
\label{constructionsect}
We fix a Hida family $\Hf$
\[
\Hf=\sum_{n=1}^{+\infty}a_n(\Hf)q^n\in \SSS^\ord(N_\Hf,\chi_\Hf,\Lambda_\Hf)
\]
primitive of tame level $N_\Hf$, tame character $\chi_\Hf$ of conductor dividing $N_\Hf\cdot p$.

We also let
\[
\Hg=\sum_{n=1}^{+\infty}a_n(\Hg)q^n\in\SSS_{\Omega_1}(M,\chi_\Hg,R_\Hg)\quad\text{and}\quad \Hh=\sum_{n=1}^{+\infty}a_n(\Hh)q^n\in\SSS_{\Omega_2}(M,\chi_\Hh,R_\Hh)
\]
be two generalized normalized $\Lambda$-adic eigenforms with $\chi_\Hf\cdot\chi_\Hg\cdot\chi_\Hh=\omega^{2a}$ for some integer $a$, where as usual $\omega$ denotes the Teichmüller character. Assume that $N_\Hf\mid M$. In the language of \cite{Hsi2021}, we are implicitly thinking about $\Hg$ and $\Hh$ as \emph{test vectors} for families of tame level dividing $M$. We also assume that $L$ contains a primitive $M$-th root of unity from now on.

\medskip
For $s\in\Z_p^\times$ and $R\in\CCC_\Lambda$ we always write $\langle s\rangle_R^{1/2}=\langle\tilde{s}\rangle_R$ where $\tilde{s}$ is the unique root of the polynomial $X^2-s\cdot\omega^{-1}(s)$ lying in $1+p\Z_p$. We also write $\langle s\rangle_R^{-1/2}=\langle s^{-1}\rangle_R^{1/2}$ (note that this does not create ambiguity).

\medskip
Let $R_{\Hf\Hg\Hh}:=\Lambda_\Hf\hat{\otimes}_{\OO_L}R_\Hg\hat{\otimes}_{\OO_L}R_\Hh$ and set
\begin{equation}
\label{bigthetacar}
    \Theta_{\Hf\Hg\Hh}:=\Theta:\Z_p^\times\to R_{\Hf\Hg\Hh}^\times\qquad\Theta(s):=\omega^{-a-1}(s)\cdot\langle s\rangle_{\Lambda_\Hf}^{1/2}\hat{\otimes}\langle s\rangle_{R_\Hg}^{-1/2}\hat{\otimes}\langle s\rangle_{R_\Hh}^{-1/2}.
\end{equation}
View $R_{\Hf\Hg\Hh}$ as $\Lambda$-algebra via $[s]\mapsto \langle s\rangle_{\Lambda_\Hf}\hat{\otimes}1\hat{\otimes}1$ for $s\in 1+p\Z_p$.

\medskip
We define a $\Theta$-twist operator on $q$-expansions given by
\begin{equation}
\label{bigthetaoperator}
|_\Theta:R_{\Hf\Hg\Hh}[\![q]\!]\to R_{\Hf\Hg\Hh}[\![q]\!]\qquad Z=\sum_{n=0}^{+\infty}a_nq^n\mapsto Z|_\Theta=\sum_{p\nmid n}\Theta(n)a_nq^n\,.
\end{equation}

Now let $\HXi:=\Hg\times(\Hh|_\Theta)$ and define
\[
\Omega^0_{\Hf\Hg\Hh}:=\{w=(x,y,z)\in\Omega_\Hf\times\Omega_\Hg\times\Omega_\Hh\mid k_x=k_y+k_z,\;k_z\geq 2\}
\]
One checks that for $w=(x,y,z)\in\Omega^0_{\Hf\Hg\Hh}$ it holds
\[
(\Hh|_\Theta)_w=\Hh_z\otimes\psi_w\in S_{k_z}(Mp^?,\chi_\Hh\omega^{2-k_z}\varepsilon_z\psi_w^2,\C_p)\,,
\]
where (for $(n,p)=1$) we set
\[
\psi_w(n)=\omega^{-a-1}(n)\cdot\varepsilon_x(n\omega^{-1}(n))^{1/2}\cdot\varepsilon_y(n\omega^{-1}(n))^{-1/2}\cdot \varepsilon_z(n\omega^{-1}(n))^{-1/2}\,.
\]
It follows that
\[
\HXi_w=\Hg_y\times(\Hh_z\otimes\psi_w)\in S_{k_x}(Mp^?,\chi^{-1}_\Hf\omega^{2-k_x}\varepsilon_x,\C_p)\,.
\]
Notice that by our definition of $\Lambda$-algebra structure on $R_{\Hf\Hg\Hh}$, for $w=(x,y,z)\in\Omega_\Hf\times\Omega_\Hg\times\Omega_\Hh$ it holds $k_w=k_x$. It follows easily that $\Omega^0_{\Hf\Hg\Hh}$ is a $(\Lambda,R_{\Hf\Hg\Hh})$-admissible set of classical integral weights.

Looking at integral classical weights specializations $w\in\Omega^0_{\Hf\Hg\Hh}\cap(\Omega_{\Hf,\Z}\times\Omega_{\Hg,\Z}\times\Omega_{\Hh,\Z})$ it is easy to deduce that, according to our definitions, it holds
\[
\HXi\in\SSS(M,\chi^{-1}_\Hf,R_{\Hf\Hg\Hh})\,.
\]
Thanks to proposition \ref{propordproj}, we can thus consider the ordinary projection \[
\HXi^{\ord}:=e(\HXi)\in\SSS^{\ord}(M,\chi^{-1}_\Hf,R_{\Hf\Hg\Hh})=\SSS^\ord(M,\chi^{-1}_\Hf,\Lambda_\Hf)\otimes_{\Lambda_\Hf} R_{\Hf\Hg\Hh}\,,
\]
where the last equality follows easily from proposition \ref{classicordR} and we emphasize (again) that the structure of $\Lambda_\Hf$-algebra on $R_{\Hf\Hg\Hh}$ is given by $a\mapsto a\hat{\otimes}1\hat{\otimes}1$ for $a\in\Lambda_\Hf$.

\bigskip
We can proceed as in \cite{Hsi2021} to define the triple product $p$-adic $L$-function. We will need an assumption on our $\Hf$.
\begin{ass}[CR]
\label{CRass}
The residual Galois representation $\bar{\V}_\Hf$ of the big Galois representation $\V_\Hf$ attached to $\Hf$ is absolutely irreducible and $p$-distinguished.
\end{ass}

\medskip
Let $\Tr_{M/N_\Hf}:\SSS^\ord(M,\chi^{-1}_\Hf,\Lambda_\Hf)\to\SSS^\ord(N_\Hf,\chi^{-1}_\Hf,\Lambda_\Hf)$ be the usual trace map.

\medskip
By the primitiveness of $\Hf$ and assumption \ref{CRass}, it follows that the so-called congruence ideal $C(\Hf)\subset\Lambda_\Hf$ of $\Hf$ is principal, generated by a non-zero element $\eta_\Hf$, called the congruence number for $\Hf$ (it is unique up to units). One can prove that $\breve{\Hf\,}$ is primitive as well and that $\Hf$ and $\breve{\Hf\,}$ have the same congruence number.

\medskip
Since $\Hf$ is primitive, we also get an idempotent operator $e_\Hf$ lying in $\T^\ord_{\mmm_\Hf}\otimes_{\Lambda_\Hf}\Frac(\Lambda_\Hf)$, where $\mmm_{\Hf}$ the maximal ideal of $\T^\ord:=\T^\ord(N_\Hf,\chi_\Hf,\Lambda_\Hf)$ corresponding to $\Hf$ and $\T^\ord_{\mmm_{\Hf}}$ is the localization of $\T^\ord$ at such maximal ideal. Morally, $e_{\Hf}$ plays the role of a projection to the $\Hf$-Hecke eigenspace. A similar discussion applies to $\breve{\Hf\,}$.

\medskip
Then we can let $e_{\breve{\Hf\,}}$ act on $\SSS^\ord(N_\Hf,\chi^{-1}_\Hf,\Lambda_\Hf)\otimes_{\Lambda_\Hf}\Frac(\Lambda_\Hf)$ and, by definition of congruence number, one has that $\eta_\Hf\cdot e_{\breve{\Hf\,}}(\Hxi)\in\SSS^\ord(N_\Hf,\chi^{-1}_\Hf,\Lambda_\Hf)$ for all $\Hxi\in\SSS^\ord(N_\Hf,\chi^{-1}_\Hf,\Lambda_\Hf)$.

\medskip 
We refer to \cite[section 3.3]{Hsi2021} and to \cite[section 3.5]{Col2020} for a more detailed discussion concerning congruence numbers and idempotents attached to primitive Hida families.

\begin{defi}
\label{padicLfunction}
With the above notation, the generalized $\Hf$-unbalanced triple product $p$-adic $L$-function $\Ls^f_p(\Hf,\Hg,\Hh)$ attached to the triple $(\Hf,\Hg,\Hh)$ is defined as
\[
\Ls^f_p(\Hf,\Hg,\Hh):=a_1\left(\eta_\Hf\cdot e_{\breve{\Hf\,}}\left(\Tr_{M/N_\Hf}\left(\HXi^\ord\right)\right)\right)\in R_{\Hf\Hg\Hh}.
\]
\end{defi}

\begin{rmk}
We view $\Ls^f_p(\Hf,\Hg,\Hh)$ as a function on $\W_{\Lambda_\Hf}(\C_p)\times\W_{R_\Hg}(\C_p)\times\W_{R_\Hh}(\C_p)$. In particular for $w=(x,y,z)\in\Omega_{\Hf\Hg\Hh}$ one gets that the evaluation of $\Ls^f_p(\Hf,\Hg,\Hh)$ at $w$ is given by
\[
\Ls^f_p(\Hf,\Hg,\Hh)(w)=\eta_{\Hf_{\!x}}\cdot a_1(e_{\breve{f}}(\Tr_{M/N_\Hf}(\HXi^\ord_w))).
\]
Recall that $(\Hh|_\Theta)_w$ is in the image of the $m=(k_x-k_y-k_z)/2$-th power of Serre's derivative operator $d=q\tfrac{d}{dq}$ acting on $p$-adic modular forms of weight $k_z$, where if $m$ is negative one defines the $m$-th power of $d$ as a $p$-adic limit. We can conclude that $(\Hh|_\Theta)_w$ is the $q$-expansion of a $p$-adic modular form of weight $k_x-k_y$ and tame level $M$. Hence, by Hida's classicality theorem for ordinary forms, we deduce that
\[
\HXi^\ord_w=e(\Hg_y\times d^m(\Hh_z\otimes\psi_w))\in S^\ord_{k_x}(Mp^t,\chi^{-1}_\Hf\omega^{2-k_x}\varepsilon_x,\C_p)
\]
where $\psi_w=\omega^{-a-1-m}\varepsilon_x^{1/2}\varepsilon_y^{-1/2}\varepsilon_z^{-1/2}$ and $t\geq 1$ depends on $w$, $\chi_{\Hg}$ and $\chi_{\Hh}$ (and it is always chosen to be large enough).
\end{rmk}

\subsection[\texorpdfstring{Evaluation of the $p$-adic $L$-function in terms of Petersson products}{Evaluation of the p-adic L-function in terms of Petersson products}]{Evaluation of the $p$-adic $L$-function in terms of Petersson products}
\begin{defi}
\label{Peterssonscalarproduct}
We set our conventions for the Petersson inner product on the spaces $S_k(N,\chi)$ of complex modular forms of level $N$ and character $\chi$ to be
\[
\langle \xi_1, \xi_2\rangle_{Pet} :=\frac{1}{\Vol(\HH/\Gamma_0(N))}\int_{\DD_0(N)}\xi_1(\tau)\overline{\xi_2(\tau)}v^k\frac{du\,dv}{v^2}
\]
for $\xi_1,\, \xi_2\in S_k(N,\chi)$ where we write $\tau=u+iv\in\HH$ (the upper half-plane) and $\DD_0(N)$ is a fundamental domain for the action of $\Gamma_0(N)$ on $\HH$.
\end{defi}

\begin{rmk}
Note that by the above definition our Petersson inner product is linear in the first variable and conjugate linear in the second variable. Moreover, it is normalized so that it does not depend on the level $N$ considered.
\end{rmk}

\begin{prop}
\label{padicLvalues}
Pick $w=(x,y,z)\in\Omega_{\Hf\Hg\Hh}$ and set
\[
C:=C_{N_\Hf,M}:=[\Gamma_0(N_\Hf):\Gamma_0(M)]=\tfrac{M}{N_\Hf}\cdot\prod_{\substack{\ell\mid M \\ \ell\nmid N_\Hf}}\left(1+\tfrac{1}{\ell}\right)\in\Z_{\geq 1}.
\]
Write
$f=\Hf_{\!x}$, $\breve{f}=(\breve{\Hf\,})_x$, $\Xi=\HXi^\ord_w\in S^\ord_{k_x}(Mp^t,\chi^{-1}_\Hf\omega^{2-k_x}\varepsilon_x,\C_p)$ to simplify the notation, so that $\breve{f}=\lambda_N(f)^{-1}\cdot w_N(f)$ as before. Assume that $t\geq1$ is large enough (in particular larger that the $p$-order of the exact level of $f$). Then the evaluation of $\Ls^f_p(\Hf,\Hg,\Hh)$ at $w$ can be described as follows, depending on two mutually exclusive cases. 
\begin{itemize}
    \item [(A)] Assume $f$ is a newform in $S_k(N_\Hf\, p^s,\chi_\Hf\omega^{2-k}\varepsilon,L)$. Then:
    \begin{equation}
    \label{computcaseC}
        \Ls^f_p(\Hf,\Hg,\Hh)(w) =\frac{\eta_{f}\cdot C\cdot p^{k(t-s)}}{a_p(\breve{f})^{t-s}}\cdot\frac{\langle\Xi,V_p^{t-s}(\breve{f})\rangle_{Pet}}{\| f\|^2_{Pet}}\,.
    \end{equation}
    \item [(B)] Assume that $f$ is the ordinary $p$-stabilization of a newform $f^\circ\in S_k(N_\Hf,\chi_\Hf^\circ,L)$ (where $\chi_\Hf^\circ$ is the $N_\Hf$-part of $\chi_\Hf$). 
    Set $f^\#:=w_{N_\Hf p}(\breve{f}^\rho)$, where $\breve{f}^\rho$ is obtained from $\breve{f}$ applying complex conjugation to the Fourier coefficients. Then
    \begin{equation}
    \label{computcaseB}
        \Ls^f_p(\Hf,\Hg,\Hh)(w)=\frac{\eta_{f}\cdot C\cdot p^{k(t-1)}}{a_p(\breve{f})^{t-1}}\cdot\frac{\langle\Xi,V_p^{t-1}(f^\#)\rangle_{Pet}}{\langle \breve{f},f^\#\rangle_{Pet}}\,.
    \end{equation}
\end{itemize}
\begin{proof}
This follows directly from \cite[proposition 4.5]{Hi1985a} (note that our conventions for the Petersson inner product differ from those of Hida, so we have to adjust the result accordingly).
\end{proof}
\end{prop}

\begin{rmk} In case $(B)$ of the above proposition (with the notation as above), assume that $t=1$ and that we can write
\[
e_{\breve{f}}(\Tr_{Mp^t/N_\Hf p^t}(\Xi))=\xi-\beta_k\chi_\Hf^\circ(p)^{-1}\cdot V_p(\xi)
\]
for some $\xi\in S_k(N_\Hf,(\chi^\circ_\Hf)^{-1})$. Then one can check that
\[
\frac{\langle\Xi,f^\#\rangle_{Pet}}{\langle \breve{f},f^\#\rangle_{Pet}} = \frac{\langle \xi,\breve{f}^\circ\rangle_{Pet}}{\langle f^\circ,f^\circ\rangle_{Pet}} =\frac{\langle\Xi,\breve{f}\rangle_{Pet}}{\langle f,f\rangle_{Pet}}.
\]
In particular, assume that $\Hg$ and $\Hh$ are classical Hida families of tame level $N_\Hf$ with $\chi_\Hf\chi_{\Hg}\chi_{\Hh}=\mathbbm{1}$ and $w=(k,l,m)\in\Omega_{\Hf\Hg\Hh}$ is a triple of classical integral weights such that $\Hg_l$ and $\Hh_m$ are ordinary $p$-stabilizations of forms $g^\circ\in S_l(N_\Hf,\chi_{\Hg}^\circ)$ and $h^\circ\in S_m(N_\Hf,\chi_{\Hh}^\circ)$ respectively. Then the hypothesis made on $\Xi$ is verified (cf. \cite[section 4.4]{BSV2020a}) and we recover the $p$-adic periods which are denoted by
$I_p(f^\circ,h^\circ,g^\circ)$ in \cite[section 1.1]{BSV2020a} and by $\Ls_p^f(f_\alpha,h_\alpha,g_\alpha)$ in \cite[section 3.1]{BSVast1}. Note that we have switched the role of $\Hg$ and $\Hh$ in our construction, compared to what happens in \cite{BSV2020a} and \cite{BSVast1}.
\end{rmk}

\subsection[\texorpdfstring{Comparison with the complex $L$-values}{Comparison with complex L-values}]{Comparison with the complex $L$-values}
In this section we compare the values of our square root triple product $L$-function with the central values of the Garret-Rankin triple product $L$-function associated to a triple of modular forms. Most of the material contained in this section is derived from \cite[section 3]{Hsi2021}.

\medskip
In this section we fix positive integers $N,M$ coprime to $p$ such that $N\mid M$. We consider a triple of cuspidal modular forms
\[
f=\sum_{n=1}^{+\infty}a_n(f)q^n,\qquad g=\sum_{n=1}^{+\infty}a_n(g)q^n,\qquad h=\sum_{n=1}^{+\infty}a_n(h)q^n
\]
with
\[
f\in S_k(Np^{e_1},\chi_f\omega^{2-k}\varepsilon_1),\; g\in S_l(Mp^{e_2},\chi_g\omega^{2-l}\varepsilon_2),\; h\in S_m(Mp^{e_3},\chi_h\omega^{2-m}\varepsilon_3),
\]
where $e_i\geq 1$ and $\varepsilon_i$ are Dirichlet characters of $p$-power order for $i=1,2,3$, while $\chi_f$ (resp. $\chi_\xi$ for $\xi\in\{g,h\}$) is a Dirichlet character defined modulo $Np$ (resp. $Mp$). 

\begin{ass}
\label{analyticass}
\begin{itemize}
    \item [(i)] $f,g,h$ are normalized eigenforms, i.e., for $\xi\in\{f,g,h\}$ it holds $a_1(\xi)=1$ and $\xi$ is an eigenform for all the Hecke operators $T_\ell$ for all primes $\ell\nmid N$ (resp. $\ell\nmid M$ if $\xi\in\{g,h\}$). We also assume that $f,g,h$ are eigenforms for the $U_p$ operator. 
    \item [(ii)] The triple $(f,g,h)$ is tamely self-dual, i.e., $\chi_f\cdot\chi_g\cdot\chi_h=\omega^{2a}$ for some integer $a$.
    \item [(iii)] The triple of weights $(k,l,m)$ is arithmetic and $f$-unbalanced, i.e., $\nu\geq 1$ for $\nu\in\{k,l,m\}$,  $k+l+m$ is even and $k\geq l+m$.
    \item [(iv)] The form $f$ is a $p$-stabilized ordinary newform, i.e., it is either the ordinary $p$-stabilization of a $p$-ordinary newform $f^\circ$ of level $N$ or an ordinary newform of level $Np^{e_1}$.
    \item [(v)] The tame level $N$ is a squarefree integer.
\end{itemize}
\end{ass}

When $f$ is the ordinary $p$-stabilization of a newform $f^\circ$ of level $N$, we write $\alpha_f$, $\beta_f$ for the roots of the Hecke polynomial at $p$ for $f^\circ$ and we always assume that $|\alpha_f|_p=1$.

\medskip
Let $r=(k+l+m)/2$ and let $\chi_\Ad$ be the adèlization of the Dirichlet character
\[
\chi:=\omega^{a-r}(\varepsilon_1\varepsilon_2\varepsilon_3)^{1/2}.
\]
Let $\pi_1=\pi_f\otimes\chi_\Ad$, $\pi_2=\pi_g$, $\pi_3=\pi_h$, where for $\xi\in\{f,g,h\}$ we denote by $\pi_\xi$ the irreducible automorphic representation of $\GL_2(\Ad)$ associated to $\xi$ as in \cite[chapter 3]{Bum1997}.

\medskip
It is well-known that there is a decomposition $\pi_\xi=\bigotimes_{\ell\leq\infty}\pi_{\xi,\ell}$ into local representations.

\medskip
Finally, let $\Pi:=\pi_1\times\pi_2\times\pi_3$ denote the corresponding automorphic representation of $\GL_2(\Ad_E)$ where $E=\Q\times\Q\times\Q$ is the split cubic étale algebra over $\Q$. Thanks to our choices one can verify that the central character of $\Pi$ is trivial, so that $\Pi$ is isomorphic to its contragradient.

\medskip
We let $L(\Pi,s)$ denote the triple product complex $L$-function attached to $\Pi$ (cf. for instance \cite{PSR1987}). It is known (cf. for instance the summary in \cite[pagg. 225-228]{Ike1992} and the references therein) that $L(\Pi,s)$ is given by a suitable Euler product converging for $\re(s)\gg 0$ and that it admits analytic continuation to an entire function with a functional equation of the form
\[
L^*(\Pi,s)=\varepsilon(\Pi,s)\cdot L^*(\Pi,1-s)
\]
Here $L^*(\Pi,s)=L(\Pi,s)\cdot L(\Pi,s)_\infty$ with
\[
L(\Pi,s)_\infty=\Gamma_\C(s+r-3/2)\cdot\Gamma_\C(s-r+k+1/2)\cdot\Gamma_\C(s+r-l-1/2)\cdot\Gamma_\C(s+r-m-1/2)
\]
and $\Gamma_\C(s)=2(2\pi)^{-s}\Gamma(s)$ ($\Gamma(\cdot)$ being Euler's gamma function). This explicit description of the archimedean $L$-factor is proven in \cite{Ike1998}.

\medskip
Moreover, $\varepsilon(\Pi,s)=\prod_{\ell\leq\infty}\varepsilon_\ell(\Pi,s)$ is an invertible function satisfying the property that $\varepsilon_\ell(\Pi,1/2)\in\{\pm 1\}$ and $\varepsilon_\ell(\Pi,1/2)=1$ for almost all $\ell$. In particular, it is known that:
\begin{itemize}
    \item [(a)] $\varepsilon_\infty(\Pi,1/2)=1$ in our case (this depends on the fact that the triple of weights $(k,l,m)$ is unbalanced);
    \item [(b)] $\varepsilon_\ell(\Pi,1/2)=1$ if $\ell\nmid pM$.
\end{itemize}

We are then led to the following further assumption.
\begin{ass}
\label{localsigns}
In what follows we assume that $\varepsilon_\ell(\Pi)=1$ for all $\ell\mid M$.
\end{ass}

\begin{defi}
If $\pi$ is an irreducible smooth representation of $\GL_2(\Q_\ell)$ for a rational prime $\ell$ and $\VV_\pi$ is a realization of $\pi$, we let $c(\pi)$ denote the smallest integer (which exists, by smoothness) such that $\VV_\pi^{\UU_1(\ell^{c(\pi)})}\neq 0$, where for all $m\in\Z_{\geq 0}$ we set
\[
\UU_1(\ell^m):=\left\{\begin{psmallmatrix} a & b\\ c & d
\end{psmallmatrix}\in\GL_2(\Z_\ell)\mid \ord_\ell(c)\geq m,\,\ord_\ell(d-1)\geq m\right\}\,.
\]
\end{defi}

\medskip
Now we connect this discussion to the triple product $p$-adic $L$-function, assuming that $f=\Hf_{\!x}$, $g=\Hg_y$, $h=\Hh_z$ are suitable specializations of families of the types considered in section \ref{constructionsect} with $w=(x,y,z)\in\Omega_{\Hf\Hg\Hh}$ so that $k_x=k, k_y=l, k_z=m$ (with $k\geq l+m$ as we have assumed before). Write $\Pi_w$ for the corresponding automorphic representation of $\GL_2(\Ad_E)$.

\medskip
Following Harris-Kudla (\cite{HK1991}) and Ichino (\cite{Ich2008}), Hsieh proved in \cite{Hsi2021} the following fact.

\begin{prop}
\label{unbalinterpolation}
Under assumptions \ref{analyticass} and \ref{localsigns}, the following formula holds:
\begin{equation}
\label{classinterp}
  (\Ls_p^f(\Hf,\Hg,\Hh)(w))^2=\chi_{f,\,p}(-1)\cdot(-1)^{k+1}\cdot\frac{L^*(\Pi_w,1/2)}{\zeta_\Q(2)^2\cdot\Omega_f^2}\cdot\mathscr{I}^{unb}_{\Pi_w,p}\cdot\left(\prod_{\ell\mid M}\mathscr{I}^*_{\Pi_w,\ell}\right)
\end{equation}
where
\begin{itemize}
    \item [(i)] \begin{equation}
        \Omega_f:=2\cdot (-2\sqrt{-1})^{k+1}\cdot\|f^\circ\|^2_{Pet}\cdot\EE_p(f,\mathrm{Ad})\cdot\eta_f^{-1}\cdot [\SL_2(\Z):\Gamma_0(N_\Hf p^{c(\pi_{f,p})} )]
    \end{equation}
    with
    \[
    f^\circ:=\begin{cases}
    f  & \text{in case (A) of prop. \ref{padicLvalues}}\\
    \text{the newform of level }N_\Hf\text{ associated to }f & \text{in case (B) of prop. \ref{padicLvalues}}    
    \end{cases}
    \]
    and
   \begin{equation}
        \EE_p(f,\mathrm{Ad})=a_p(f)^{-c(\pi_{f,p})}\cdot p^{c(\pi_{f,p})(k/2-1)}\cdot\varepsilon(\pi_{f,p},1/2)\cdot\sigma_f,
    \end{equation}
    where
    \[
    \sigma_f:=\begin{cases}
    1 & \text{in case (A) of prop. \ref{padicLvalues}}\\
    \left(1-\tfrac{\beta_f}{\alpha_f}\right)\left(1-\tfrac{\beta_f}{p\alpha_f}\right) & \text{in case (B) of prop. \ref{padicLvalues}, equiv. if }c(\pi_{f,p})=0    
    \end{cases}
    \]
    \item [(ii)] $\mathscr{I}^{unb}_{\Pi_w,p}$ is the normalized local zeta integral defined as \cite[equation 3.28]{Hsi2021};
    \item [(iii)] $\mathscr{I}^*_{\Pi_w,\ell}$ is the normalized local zeta integral defined as in \cite[equation 3.29]{Hsi2021};
    \item [(iv)] $\zeta_\Q(\cdot)=\pi^{-1}\zeta(\cdot)$ where $\zeta(\cdot)$ is the usual Riemann zeta function, so that $\zeta_\Q(2)=\pi/6$;
    \item [(v)] $\chi_{f,\,p}$ denotes the $p$-part of the character $\chi_f$.
\end{itemize}
\begin{proof}
This is essentially a restatement of proposition 3.10 and corollary 3.13 in \cite{Hsi2021}. Note that our normalization for the Petersson inner product is different from Hsieh's. This explains the appearance of the factor $\zeta_\Q(2)^2$ in our formula and the slight changes in the definition of the period $\Omega_f$.
\end{proof}
\end{prop}

\begin{rmk}
One can compute directly  that, if we are in case (B) of proposition \ref{padicLvalues}, it holds that
\[
\|f^\circ\|^2_{Pet}\cdot\sigma_f\cdot\frac{(-1)^k\cdot\alpha_f\cdot\chi_\Hf^\circ(p)^{-1}}{\lambda_N(f)\cdot p^{k/2}\cdot(1+1/p)}=\langle \breve{f}, f^\#\rangle_{Pet}.
\]
We refer \cite[proposition 5.4.1]{Col2020} for a very similar computation, where the form denoted $h^{\musNatural{}}$ there should be thought as a constant multiple of our $f^\#$. 
This explains the appearance of the factor $\sigma_f$ and allows an even more direct comparison (in the $\Hf$-unbalanced region) between the formula given by equation \ref{classinterp} and the formulas appearing in the statement of proposition \ref{padicLvalues}.
\end{rmk}

\begin{rmk}
\label{canonicalcongrnumber}
In the recent preprint \cite{Mak2023}, Maksoud has defined a canonical congruence number attached to a Hida family $\Hf$ (satisfying assumption \ref{CRass}) with precise interpolation properties. Such congruence number is interpreted as a \emph{weight variable} adjoint $p$-adic $L$-function of $\Hf$ and denoted $L_p(\ad\Hf)$. Fix $x\in\Omega_\Hf$ and write $f=\Hf_{\!x}$. Assume that the exact level of $f$ is $N_\Hf p^r$ and that $x$ lies over $(k,\varepsilon)\in\W_\Lambda$. Write $N_\Hf p^{r_0}$ for the exact level of the newform $f^\circ$ attached to $f$ as in the statement of proposition \ref{unbalinterpolation} (i.e., $r_0=c(\pi_{f,p})$). Then \cite[theorem 2.3.7]{Mak2023} claims that:
\begin{align*}
 L_p & (\ad\Hf)(x)=\\
&=p^{r-1}\cdot a_p(f)^{r}\cdot\lambda_{N_\Hf p^{r_0}}(f^\circ)\cdot\tilde{\sigma}_f\cdot (-2i)^{k+1}\cdot\frac{\zeta_\Q(2)\cdot\| f^\circ\|^2_{Pet}}{\Omega_{f^\circ}^+\cdot\Omega_{f^\circ}^-}\cdot[\SL_2(\Z):\Gamma_0(Np^{r_0})]\,,   
\end{align*}
where:
\begin{itemize}
    \item [(i)] \[
    \tilde{\sigma}_f:=\begin{cases}
    1 & \text{in case (A) of prop. \ref{padicLvalues}}\\
    \alpha_f\cdot (p-1)\cdot \left(1-\tfrac{\beta_f}{\alpha_f}\right)\left(1-\tfrac{\beta_f}{p\alpha_f}\right) & \text{in case (B) of prop. \ref{padicLvalues}}   
    \end{cases}\,;
    \]
    \item [(ii)] the periods $\Omega_{f^\circ}^\pm\in\C^\times$ are the usual complex periods attached to eigenforms, whose definition is recalled in \cite[section 2.1]{Mak2023}.
\end{itemize}

\medskip
In what follows we fix the choice $\eta_\Hf=L_p(\ad\Hf)\in\Lambda_\Hf$ for the congruence number of $\Hf$. We leave to the interested reader the task of applying the interpolation formula for $L_p(\ad\Hf)$ in order to make the formulas involving $\eta_\Hf$ more explicit.
\end{rmk}

\section{Families of theta series of infinite \texorpdfstring{$p$}{p}-slope}
\label{CMfamilies}
\subsection{Setup for the interpolation}
\label{setup}
We fix an odd prime $p$ and we let $K$ be an imaginary quadratic field where $p$ is inert. Denote by $N_{K/\Q}$ the norm morphism on fractional ideals in $K$. Let $-d_K$ be the discriminant of $K$ (so that $p\nmid d_K$) and let $\varepsilon_K$ denote the central character of $K$, i.e., more explicitly
\[
\varepsilon_K(n)=\left(\frac{-d_K}{n}\right)\qquad\text{if}\quad (n,d_K)=1
\]
where $\left(\frac{\;\cdot\;}{\;\cdot\;}\right)$ denotes the Jacobi symbol.

\begin{defi}
For $\aid\subset\OO_K$ an integral ideal in $\OO_K$, we let $I_K(\aid)$ denote the group of fractional ideals of $K$ prime to $\aid$ and we set
\[
P_K(\aid):=\{(\alpha)\in I_K(\aid)\mid \alpha\equiv 1\mod^\times\aid\}\, ,\qquad Cl_K(\aid):=I_K(\aid)/P_K(\aid).
\]
The group $Cl_K(\aid)$ is the so-called ray class group modulo $\aid$.   
\end{defi}

\begin{rmk}
It is well-known that $Cl_K(\aid)$ is a finite group.
\end{rmk}

We fix a finite order character $\eta: G_K\to\bar{\Q}^\times$ with conductor $\cc$ (a non-trivial proper integral ideal in $\OO_K$). Via class field theory we will freely view $\eta$ as a \emph{ray class character} $\eta:Cl_K(\cc)\to\bar{\Q}^\times$ or a finite order character $\eta:\Ad_K^\times/K^\times\to \bar{\Q}^\times$ (note the slight abuse of notation here). Moreover, we assume that $\eta$ is not the restriction of a character of $G_\Q$.

\medskip
Denote by $\eta_{|\Q}$ the Dirichlet character defined modulo $N_{K/\Q}(\cc)$ and given by
\[
\eta_{|\Q}(n):=\eta((n))\qquad\text{for}\quad (n,N_{K/\Q}(\cc))=1
\]
It is then a classical theorem of Hecke and Shimura (cf. \cite{Miy1989}, theorem 4.8.2) that the $q$-expansion (where as usual $q=\exp(2\pi i \tau)$ for $\tau\in\HH$)
\begin{equation}
\label{thetaseries}
    g(\tau):=\theta_\eta(\tau):=\sum_{(\aid,\cc)=1}\eta(\aid)q^{N_{K/\Q}(\aid)}
\end{equation}
defines a cuspidal modular form of weight $1$ (the theta series attached to the character $\eta$). Here the sum runs over the integral ideals in $\OO_K$ prime to $\cc$. 

More precisely, $g\in S_1(d_K\cdot N_{K/\Q}(\cc), \varepsilon_K\cdot\eta_{|\Q})$ and since we assume that $\eta$ is of exact conductor $\cc$, $g$ is also a newform of level $d_K\cdot N_{K/\Q}(\cc)$. From now on, we set $N_g:=d_K\cdot N_{K/\Q}(\cc)$ and $\chi_g:=\varepsilon_K\cdot\eta_{|\Q}$.

\medskip
The Fourier coefficients of $g$ generate a finite extension of $\Q$. We can thus view $g$ as a modular form whose $q$-expansion at $\infty$ has coefficients in a finite extension $L$ of $\Q_p$ (via the embedding $\iota_p$), i.e., $g\in S_1(N_g,\chi_g,L)$. As in the previous sections, we assume that $L$ is large enough. In particular, here we assume that $L$ contains the completion of $K$ inside $\C_p$ (which we will denote by $K_p$ with ring of integers $\OO_{K,p}$).

\medskip
We would like to find a $p$-adic family of modular forms - all with complex multiplication by K - of varying weights (in the sense of Hida-Coleman) having $g$ (or a slight modification of $g$) as a specialization in weight $1$. We will see that this can actually be done explicitly.

\begin{rmk}
Since the fixed prime $p$ is inert in $K$, $p^r\OO_K\mid\cc$ if and only if $p^{2r}\mid N_{K/\Q}(\cc)$. Hence we should distinguish two cases:
\begin{itemize}
\item[(a)] $(p\OO_K,\cc)=1$, or equivalently $p\nmid N_g$
\item[(b)] $\ord_p(N_g)=2r$ for some $r\in\Z_{\geq 1}$
\end{itemize}
In both cases it holds that $a_p(g)=0$, or equivalently that $T_p(g)=0$ in case (a) (resp. $U_p(g)=0$ in case (b)). This is usually described as $g$ having infinite $p$-slope.
\end{rmk}

\begin{rmk}
While case (a) can be reinterpreted in the realm of Hida theory (as in this case $g$ admits one or two ordinary $p$-stabilizations), case (b) is instead more genuinely a problem in \emph{infinite slope}. This dichotomy is also reflected in the fact that the local component at $p$ of the automorphic representation associated with $g$ is a principal series in case (a) and a supercuspidal representation in case (b).    
\end{rmk}

\begin{ass}
From now on in this section we will always assume that $p\OO_K\mid\cc$ and we will write $\cc=\cc_0\cdot p^r\OO_K$ with $\cc_0$ coprime to $p\OO_K$ and $r\geq 1$.
    
\end{ass}
 
\begin{rmk}
When $p$ splits in $K$ one can explicitly write down families of theta series, specializing to ($p$-stabilizations) of  modular forms of the shape described in \eqref{thetaseries}. See, for instance, \cite[section 4.2]{BDV22} for a discussion about this construction, which - again - is well-understood within Hida theory.
\end{rmk}

\medskip
In what follows, we try to adapt such construction to our setting. Notice that $K_p/\Q_p$ is the unique degree two unramified extension of $\Q_p$ inside our fixed algebraic closure $\bar{\Q}_p$, so we will identify $K_p=\Q_{p^2}$ (with ring of integers $\Z_{p^2}$). Moreover we have a decomposition
\[
\Z_{p^2}^\times=\mu_{p^2-1}\times (1+p\Z_{p^2})
\]
induced by the Teichm\"{u}ller lift. Note that $1+p\Z_{p^2}$ does not contain $p$-power roots of unity.

\medskip
Let $G_p$ be the subgroup of the idèlic class group $C_K:=\Ad_K^{\times}/K^\times$ over $K$ defined by
\[
G_p:=K^\times\!\cdot(\C^\times\cdot\mu_{p^2-1}\cdot\prod_{\lid\neq p\OO_K}\OO_\lid^\times)/K^\times\,.
\]
Set moreover $I_{K,\infty}:=K^\times\!\cdot(\C^\times\cdot\prod_{\lid}\OO_\lid^\times)/K^\times$ and let $\Pic(\OO_K)$ denote the classical ideal class group of $K$.

The snake lemma applied to the following diagram with exact rows
\begin{center}
    \begin{tikzcd}
        & 0\arrow[r] & G_p\arrow[r]\arrow[d, hookrightarrow] & C_K\arrow[r]\arrow[d, equals] & C_K/G_p\arrow[r]\arrow[d, twoheadrightarrow] & 0\\
        &0\arrow[r] &I_{K,\infty}\arrow[r] &C_K\arrow[r] &\Pic(\OO_K)\arrow[r] &0
    \end{tikzcd}
\end{center}
identifies $1+p\Z_{p^2}\cong\Ker(C_K/G_p\twoheadrightarrow\Pic(\OO_K)$. We can thus consider the diagram
\begin{center}
\begin{tikzcd}
& 1\arrow[r] & 1+p\Z_{p^2} \arrow[r]\arrow[d, hookrightarrow, "\iota"] & C_K/G_p \arrow[r]\arrow[dl, dashed] &\Pic(\OO_K) \arrow[r] & 1\\
& & \bar{\Q}_p^\times
\end{tikzcd}
\end{center}
where the horizontal row is an exact sequence of abelian groups, $\iota$ is given by $\iota(u)=u^{-1}$ and the dashed arrow is any (continuous) extension of $\iota$ to the quotient $\Ad_K^\times/G_p$, obtained using the divisibility of $\bar{\Q}_p^{\times}$. Finally we let $\lambda^{(p)}$ to be the following composition:
\begin{center}
    \begin{tikzcd}
    & \lambda^{(p)}:\Ad_K^\times/K^\times\arrow[r, twoheadrightarrow] & \Ad_K^\times/G_p\arrow[r, dashed] & \bar{\Q}_p^\times.
    \end{tikzcd}
\end{center}

We associate to $\lambda^{(p)}$ an algebraic Hecke character of $K$ of $\infty$-type $(1,0)$ as follows:
\[
 \lambda^{(\infty)}:\Ad_K^\times/K^\times \to\C^\times\qquad 
    x=[(x_\nu)_\nu] \mapsto\left(\iota_\infty\circ\iota_p^{-1}(\lambda^{(p)}(x)\cdot x_p)\right)\cdot x_\infty^{-1}\,.
\]
Finally, writing $\lambda^{(\infty)}=\otimes_v\lambda^{(\infty)}_v$ one gets a character at the level of fractional ideals
\[
 \lambda: I_K(p\OO_K) \to\bar{\Q}^\times\qquad \aid \mapsto \prod_{\lid\mid\aid}\lambda^{(\infty)}_\lid(\varpi_\lid)^{\ord_\lid(\aid)}\,,
\]
where $\varpi_\lid$ is a uniformizer at $\lid$. One can verify that $\lambda((\alpha))=\alpha$ whenever $\alpha\equiv 1\mod^\times p\OO_K$.

\begin{defi}
In the above setting, we will say that $\lambda^{(p)}$ is the \textbf{$p$-adic avatar} of $\lambda$ and that $\lambda^{(\infty)}$ is the \textbf{complex avatar} of $\lambda$. 
\end{defi}

\begin{rmk}
We will also look at $\lambda^{(p)}$ as a $p$-adic Galois character $\lambda^{(p)}:G_K\to\bar{\Q}_p^\times$ via global class field theory.  
\end{rmk}

\medskip
Up to enlarging $L$, we can assume that $\lambda(\aid)\in L$ for all $\aid\in I_K(p\OO_K)$ and $\eta(\aid)\in L$ for all $\aid\in I_K(\cc)$. 

\begin{defi}
We let $\langle\,\cdot\,\rangle :\OO_L^\times\to\OO_L^\times$ to be the projection onto the free units (note that now $\OO_L^\times$ might contain $p$-power roots of unity). By slight abuse of notation we will write $\langle \lambda(\aid)\rangle$ to denote $\iota_p^{-1}(\langle \iota_p(\lambda(\aid))\rangle)$ (notice that this makes sense).    
\end{defi}

\begin{defi}
\label{defietakgk}
For $k\in\Z_{\geq 1}$, let $\eta_k:I_K(\cc) \to\bar{\Q}^\times$ be the character $\aid \mapsto \eta(\aid)\cdot \langle \lambda(\aid)\rangle^{k-1}$,
so that
\[
g_k:=\sum_{(\aid,\cc)=1}\eta_k(\aid)q^{N_{K/\Q}(\aid)}\in S_k(N_\Hg,\chi_k)
\]            
where $N_\Hg=N_g$ and $\chi_k=\chi_g\cdot\omega^{1-k}=\chi_\Hg\cdot\omega^{2-k}$ where $\omega$ is the Teichm\"{u}ller character and clearly $\chi_\Hg=\chi_g\cdot\omega^{-1}$. We will also write $N_\Hg^\circ:=N_\Hg/p^{2r}$ in the sequel.    
\end{defi}

\begin{rmk}
\label{notprimchar}
Note that, since $p$ is inert in $K$, the $p$-part of the conductor of $\chi_k$ is at most $p^r$ for all $k\geq 1$, so that $\chi_k$ will never be $p$-primitive as a Dirichlet character modulo $N_\Hg$. This is a typical feature for newforms of infinite $p$-slope and level divisible by $p$. It is well-known, on the other hand, that if the $p$-order of $N$ and of $\mathrm{cond}(\chi)$ of a normalized newform $f\in S_k(N,\chi)$ coincide, then $a_p(f)$ must have euclidean absolute value $p^{(k-1)/2}$ (cf. theorem 4.6.17 of \cite{Miy1989}).
\end{rmk}

\begin{rmk}
Recall the (unique) continuous $\Z_p$-action on $U_1:=\{z\in\C_p\mid |z-1|_p<1\}$ extending the natural structure of $U_1$ as a multiplicative abelian group, namely
\[
z^s:=\sum_{n=0}^\infty\binom{s}{n}(z-1)^n\qquad z\in U_1,\, s\in\Z_p.
\]
We thus view $U_1$ as a topological $\Z_p$-module. One can show that $\mu_{p^\infty}(\C_p)$ (i.e., the subgroup of roots of unity of $p$-power order) is dense inside $U_1$. It follows that the natural action of $G_{\Q_p}$ on $U_1$ given by the $p$-adic cyclotomic character $\cyc^{(p)}:G_{\Q_p}\to\Z_p^\times$ is compatible with the action of $\Z_p^\times$, in the sense that $\sigma(z)=z^{\cyc^{(p)}(\sigma)}$ for $z\in U_1,\,\sigma\in G_{\Q_p}$.
\end{rmk}

\begin{defi}
    We define $W_K$ to be the smallest closed $\Z_p$-submodule of $U_1$ containing $\langle\lambda(\aid)\rangle$ for all $\aid\in I_K(p\OO_K)$.
\end{defi}

\begin{rmk}
Note that the notation $W_K$ makes sense, since different choices for $\lambda$ (i.e., different choices for the dashed arrow in the diagram above) differ by a finite order character, so that $W_K$ only depends on $K$ and not on $\lambda$.
\end{rmk}

\begin{lemma}
\label{Wlambda}
$W_K$ is a free $\Z_p$-module of rank $2$. If $a\in\Z_{\geq 0}$ is such that $p^a=\#(Cl_K(p\OO_K)\otimes\Z_p)$, then $w^{p^a}\in 1+p\Z_{p^2}$ for all $w\in W_K$. In particular, if $p\nmid\#(\Pic(\OO_K))$, we have $W_K=1+p\Z_{p^2}$.
\begin{proof}
Let $m=\#Cl_K(p\OO_K)$. Since $\lambda((\alpha))=\alpha$ for all $\alpha\equiv 1\mod^\times p\OO_K$ we deduce that $\langle\lambda(\aid^m)\rangle\in 1+p\Z_{p^2}$ for all $\aid\in I_K(p\OO_K)$, whence $W_K^{(m)}=\{w^m\mid w\in W_K\}\subseteq 1+p\Z_{p^2}$.

Raising to the $m/p^a$-th power is an automorphism of $W_K$ as $\Z_p$-module, hence $W_K^{(p^a)}=\{w^{p^a}\mid w\in W_K\}\subseteq 1+p\Z_{p^2}$. Finally, it is also clear that $1+p\Z_{p^2}\subseteq W_K$, which proves the statment concerning the rank of $W_K$.
\end{proof}
\end{lemma}

\begin{rmk}
\label{Zp2remark}
Denote by $\langle\lambda\rangle: G_K\twoheadrightarrow W_K$ the corresponding Galois character (given by the composition $\langle\,\cdot\,\rangle\circ\lambda^{(p)}$) and let $K_\infty$ denote the (unique) $\Z_p^2$-extension of $K$. It follows from the construction that $\langle\lambda\rangle$ factors through $\Gamma_\infty:=\Gal(K_\infty/K)$, inducing an isomorphism $\Gamma_\infty\cong W_K$. We will consider $W_K$ as a $G_\Q$-module via this isomorphism (and the $G_\Q$-action on $\Gamma_\infty$ by conjugation). In particular we have $\Gamma_\infty=\Gamma^+\times\Gamma^-$ where 
\begin{itemize}
    \item [(i)] $\Gamma^+$ is the Galois group of the cyclotomic $\Z_p$-extension of $K$, denoted by $K_\infty^+$, where complex conjugation acts as the identity;
    \item [(ii)] $\Gamma^-$ is the Galois group of the anticyclotomic $\Z_p$-extension of $K$, denoted by $K_\infty^-$, where complex conjugation acts as taking the inverse.
\end{itemize}
We will write $W_K=W_K^+\times W_K^-$ for the corresponding decomposition of $W_K$.
\end{rmk}

\subsection{Families of theta series as generalized \texorpdfstring{$\Lambda$}{lambda}-adic eigenforms} 
\label{alahida}
It is possible to realize $p$-adic families of theta series of infinite $p$-slope considered above in the following way, as suggested in Hida's blue book \cite[pagg. 236-237]{Hi1993}.

\begin{defi}
We define the $\Lambda$-algebras $\Lambda_\mathrm{Hida}:=\OO_L[\![W_K]\!]$ and $\OO_\mathrm{Hida}:=\Lambda_\mathrm{Hida}[1/p]$, with $\Lambda$-algebra structure induced by the natural inclusion $1+p\Z_p\subset W_K$.
\end{defi}

\begin{defi}
\label{defGHida}
We define
\[
\Ghid:=\Sum_{(\aid,\cc)=1}\frac{\eta(\aid)}{\langle\lambda(\aid)\rangle}[\langle\lambda(\aid)\rangle]\cdot q^{N_{K/\Q}(\aid)}\in\Lambda_\mathrm{Hida}[\![q]\!]\,,
\]
\end{defi}
where recall that $[\,\cdot\,]$ denotes group elements in $W_K$.

\medskip
Let $w:\Lambda_\mathrm{Hida}\to\C_p$ be a continuous $\OO_L$-algebra homomorphism. Assume that there exists integers $a_w\geq 1$ and $k_w\geq 1$ such that $w$ sends group elements in $[u]\in 1+p^{a_w}\Z_{p^2}\subseteq W_K$ to $u^{k_w}\in\C_p$. Then
\[
\eta_w:I_K(\cc)\to\C_p^\times\qquad\aid\mapsto \frac{\eta(\aid)}{\langle\lambda(\aid)\rangle}\cdot w([\langle\lambda(\aid)\rangle])
\]
is a primitive Hecke character of infinity type $(k_w-1,0)$ with conductor $p^{e(w,\eta)}\cc$ for a suitable integer $e(w,\eta)\geq 0$ (depending on $a_w$ and the $p$-part of $\eta$), so that
\begin{equation}
\label{specHid}
    \Ghid(w):=\Sum_{(\aid,\cc)=1}\eta_w(\aid)\cdot q^{N_{K/\Q}(\aid)}\in S_{k_w}(N_w,\chi_w,\OO_L[w])\,,
\end{equation}
where
\begin{itemize}
    \item [(i)] $N_w=d_K\cdot N_{K/\Q}(\cc)\cdot p^{2e(w,\eta)}$
    \item [(ii)] $\chi_w=\varepsilon_K\cdot\eta_{|\Q}\cdot\omega^{1-k}\cdot\varepsilon_w=\chi_\Hg\cdot\omega^{2-k}\cdot\varepsilon_w$, where $\varepsilon_w$ is an explicit character valued in $\mu_{p^{\infty}}(\C_p)$, depending on $w$.
    \item [(iii)] $\OO_L[w]$ is the finite extension of $\OO_L$ generated by the values of $w$ (one can assume that it is a cyclotomic extension of $\OO_L$ generated by a $p$-power root of unity).
\end{itemize}
 
When $w$ acts on group elements $[u]\in W_K$ as $w([u])=u^k$ for some $k\geq 1$, we recover the specialisations $\Ghid(w)=g_k$.

\medskip
One clearly has a big Hecke character
\begin{equation}
    \pmb{\eta}_\mathrm{Hida}: I_K(\cc)\to\OO_\mathrm{Hida}^\times\qquad \aid\mapsto \frac{\eta(\aid)}{\langle\lambda(\aid)\rangle}\cdot[\langle\lambda(\aid)\rangle]
\end{equation}

with associated Galois character $\pmb{\eta}_\mathrm{Hida}: G_K\to\OO_\mathrm{Hida}^\times$. Note that, by construction, $\pmb{\eta}_\mathrm{Hida}$ factors through the Galois group of the ray class field modulo $\cc_0 p^\infty$ over $K$.

\begin{defi}
We set $\V_\Ghid:=\Ind_K^\Q\,\pmb{\eta}_\mathrm{Hida}$ and we call it the big Galois representation associated with the family $\Ghid$.
\end{defi}

\begin{rmk}
By construction, it follows that for any $w$ as above, the $2$-dimensional (over $L[w]=\Frac(\OO_L[w])$) $G_\Q$-representation obtained as
\[
\V_\Ghid(w):=\V_\Ghid\otimes_{\OO_\mathrm{Hida},w}L[w]
\]
is the dual of the Deligne representation attached to the specialization $\Ghid(w)$.
\end{rmk}

Now we are ready to prove that the families of the form $\Ghid$ fit in the framework of generalized $\Lambda$-adic modular forms, as defined in section \ref{generalized}.
\begin{lemma}
The family $\Ghid$ constructed as in equation \eqref{defGHida} satisfies (with the notation introduced in section \ref{generalized} and above)
\[
\Ghid\in\SSS_{\Omega_\mathrm{Hida}}(N_\Hg^\circ,\chi_\Hg,\Lambda_\mathrm{Hida})\,,
\]
where $\Omega_\mathrm{Hida}:=\Omega_{\Ghid,\Z}$. Moreover, $\Ghid$ is a generalized $\Lambda$-adic eigenform, lying in the kernel of $U_p$.
\begin{proof}
As far as $\Ghid$ is concerned, it is enough to check that $\Omega_{\Ghid,\Z}$ is $(\Lambda,\Lambda_\mathrm{Hida})$-admissible. Condition (i) of definition \ref{decentsetwt} is clearly satisfied. For condition (ii), for every $k\geq 2$ let $w_k:\Lambda_{\mathrm{Hida}}\to\C_p$ denote the weight uniquely determined by the assignment $w_k([u])=u^k-1$ on group elements. We know that $w_k\in\Omega_{\Ghid,\Z}$ and we claim that $I:=\bigcap_{k\geq 2}\Ker(w_k)=(0)$. Since $\varpi_L\notin\Ker(w_k)$ for every $k\geq 2$, one can prove the assertion working in $\Lambda_{\mathrm{Hida}}[1/p]$, where it is easy to show that $\bigcap_{k=2}^m\Ker(w_k)[1/p]=\prod_{k=2}^m\Ker(w_k)[1/p]$ for all $m\geq 2$. Using that $\Lambda_{\mathrm{Hida}}[1/p]$ is a UFD (since $\Lambda_{\mathrm{Hida}}$ is such), one concludes that indeed it must be $I=(0)$.

\end{proof}
\end{lemma}

\begin{rmk}
The families of the form $\Ghid$ are examples of $\Lambda$-adic forms admitting classical specializations also for arithmetic weights $w$ with $k_w=1$.  
\end{rmk}


\section{Factorization of triple product \texorpdfstring{$p$}{p}-adic \texorpdfstring{$L$}{L}-functions}
\label{factorizationsection}
\subsection{Remarks on the relevant complex \texorpdfstring{$L$}{L}-functions}
\label{complexrmkfact}
In this section we recollect some facts concerning Hecke $L$-functions and Rankin-Selberg convolution that will be needed in the sequel.

\medskip
Fix $K/\Q$ a quadratic imaginary field and let $\chi_\C:\Ad_K^\times/K^\times\to\C^\times$ be an algebraic Hecke character of $\infty$-type $(a,b)$. Let $|\cdot|_{\Ad_K}$ denote the adèlic norm. Then $\chi_\C=\chi_0\cdot |\,\cdot\,|_{\Ad_K}^{(a+b)/2}$ is a unitary Hecke character (i.e. taking values in $\{z\in\C^\times\mid |z|=1\}$) and the completed $L$-function $L^*(\chi_0,s)$ attached to $\chi_0$ has meromorphic continuation and functional equation with center $s=1/2$ (cf. Tate's thesis). Note that $L^*(\chi_0,s)$ is actually an entire function if $\chi_0$ is not of the form $\chi_0=\nu\circ N_{K/\Q}$ for some Dirichlet character $\nu$.

\medskip
As explained in \cite[theorem 11.3 and proposition 12.1]{JL1970}, one can attach to $\chi_\C$ an automorphic representation $\pi(\chi)$ of $\GL_2(\Ad_\Q)$ such that $L^*(\pi(\chi),s)=L^*(\chi_0,s)$. Note that if $b=0$ and $a\geq 0$, then $\pi(\chi)$ is the automorphic representation attached to the theta series $\theta_\chi$ and $L^*(\theta_\chi,s)=L^*(\pi(\chi),s+a/2)$.

\medskip
Given two automorphic representations $\pi_1$ and $\pi_2$ of $\GL_2(\Ad)$ with central characters $\omega_1$ and $\omega_2$, one can construct - via the so-called Rankin-Selberg method - an $L$-function $L^*(\pi_1\times\pi_2,s)$, prove its meromorphic continuation and functional equation of the form
\[
L^*(\pi_1\times\pi_2,s)=\varepsilon(\pi_1\times\pi_2,s)\cdot L^*(\Tilde{\pi}_1\times\Tilde{\pi}_2,1-s)
\]
where $\Tilde{\pi}$ denotes the contragradient representation of $\pi$. The poles of $L^*(\pi_1\times\pi_2,s)$ are those of $L(\omega_1\omega_2,2s-1)$. Moreover, the $\varepsilon$-factor $\varepsilon(\pi_1\times\pi_2,s)$ is an invertible function. 

\medskip
We refer to the standard reference \cite{Jac1972} for this construction and for the definition of the local $L$-factors and $\varepsilon$-factors of such $L$-functions. The local theory is also nicely summarized in \cite[section 1]{GJ1978}). For the definition of the local $\varepsilon$-factors we always use the standard additive character of the corresponding local field and the self-dual Haar measure with respect to the standard character.

\medskip
Starting from two cuspidal eigenforms $f\in S_k(N_f,\chi_f)$ and $g\in S_l(N_g,\chi_g)$, one can also define the $L$-function $L(f\times g,s)$ more classically via an Euler product expansion (cf. \cite[section 7]{Ka2004}). If $f$ and $g$ are newforms and $k\geq l$, it holds
\[
L(f\times g,s)\cdot\Gamma_\C(s)\cdot\Gamma_\C(s-l+1)=L^*(\pi_f\times\pi_g,s-\tfrac{k+l-2}{2}).
\]
Finally, if $f\in S_k(N_f,\chi_f)$ and $\psi:\Ad_K^\times/K^\times\to\C^\times$ is an algebraic Hecke character of $\infty$-type $(a,b)$, we set
\[
L^*(f/K,\psi,s):=L^*(\pi_f\times\pi(\psi),s-\tfrac{k-1+a+b}{2}).
\]
One can write $L^*(f/K,\psi,s)=L_\infty(f/K,\psi,s)\cdot L(f/K,\psi,s)$ with archimedean $L$-factor given by
\[
L_\infty(f/K,\psi,s)=\Gamma_\C(s-\min\{a,b\})\cdot\Gamma_\C(s-\min\{k-1,|a-b|\}-\min\{a,b\}).
\]

\medskip
Assume now that we are given $f\in S_k(N_f,\chi_f)$ a Hecke eigenform and two Hecke characters $\psi_1,\psi_2$ of $K$ of $\infty$-type $(l-1,0)$ and $(m-1,0)$ respectively (here $l\geq 1,m\geq 1$), which are not induced by Dirichlet characters. Then $g=\theta_{\psi_1}$ and $h=\theta_{\psi_2}$ are cuspidal newforms, say $g\in S_l(N_l,\chi_g)$ and $h\in S_m(N_h,\chi_h)$. Assume that $\chi_f\cdot\chi_g\cdot\chi_h=1$ and consider the Garret-Rankin triple product $L$-function
\[
L^*(f\times g\times h,s)=L^*(\pi_f\times\pi(\psi_1)\times\pi(\psi_2),s-\tfrac{k+l+m-3}{2})\,.
\]

If one looks at the corresponding $\ell$-adic Galois representations for $\ell$ any rational prime, one easily deduces the following decomposition
\begin{align}\nonumber
    \label{Galoisfactor1}
    V_\ell(f)\otimes V_\ell(g)\otimes V_\ell(h)&\cong V_\ell(f)\otimes \left(\Ind_K^\Q\psi_1\psi_2\oplus\Ind_K^\Q\psi_1\psi_2^\sigma\right)\cong \\
    &\cong\left(V_\ell(f)\otimes\Ind_K^\Q\psi_1\psi_2\right)\oplus\left(V_\ell(f)\otimes\Ind_K^\Q\psi_1\psi_2^\sigma\right)\,.
\end{align}

For the sake of precision, here $V_\ell(\xi)$ denotes the dual of the Deligne representation attached to $\xi$ and we look at $\psi_1$ and $\psi_2$ as Galois characters attached to the $\ell$-adic avatars of $\psi_1$ and $\psi_2$ via class field theory.

\medskip
The decomposition \eqref{Galoisfactor1} corresponds to the following factorization of $L$-functions
\begin{equation}
\label{Lfctsfactor1}
    L^*(f\times g\times h,s)=L^*(f/K,\psi_1\psi_2,s)\cdot L^*(f/K,\psi_1\psi_2^\sigma,s)\,.
\end{equation}

\subsection{Study of the big Galois representations}
As usual, we let $L$ denote a (large enough) finite extension of $\Q_p$, containing all the needed coefficients.
\begin{set}
\label{tripleprodass}
We work in the following setting (cf. section \ref{factorizationsectionintro}).
    \begin{itemize}
    \item [(i)] We fix $\Hf\in\SSS^{ord}(N_\Hf,\mathbbm{1},\Lambda_{\Hf})$ a primitive Hida family with trivial tame character, squarefree tame level $N_\Hf$ and coefficients in $\Lambda_\Hf$ (a ring in $\CCC_\Lambda$, which is also finite flat over $\Lambda=\OO_L[\![1+p\Z_p]\!]$), satisfying assumption \ref{CRass}.
    \item [(ii)] We let $K/\Q$ denote a quadratic imaginary field of odd discriminant $-d_K$ (i.e. we have $d_K\equiv 3\mod 4$) coprime to $pN_\Hf$ such that the fixed odd prime $p$ is inert in $K$ and does not divide the class number of $K$. Writing $N_\Hf=N_\Hf^+\cdot N_\Hf^-$ where $N_\Hf^+$ is the product of prime factors of $N_\Hf$ which are split in $K$, we assume that $N_\Hf^-$ is the product of an \emph{odd} number of prime factors (\emph{Heegner hypothesis}).
    \item [(iii)] We fix two ray class characters $\eta_1$ and $\eta_2$ of $G_K$, both of conductor $c p^r\OO_K$ with $c$ a positive integer with $(c,pN_\Hf)=1$ and $r\geq 1$. We then let $\Hg$ and (respectively) $\Hh$ denote the generalized $\Lambda$-adic eigenforms attached to $\eta_1$ and (respectively) $\eta_2$ via the construction explained in section \ref{alahida}.
    \item [(iv)] We assume that the central characters of $\eta_1$ and $\eta_2$ are inverse to each other (\emph{self-duality condition}).
    \item [(v)] We assume that the prime divisors of the integer $c$ are all split in $K$.
\end{itemize}
\end{set}

\medskip
\begin{rmk}
\label{cyclolambda}
Let $\langle\cyc\rangle:G_\Q\to 1+p\Z_p$ be the character $g\mapsto\cyc(g)\cdot\omega(\cyc(g))^{-1}$, where $\cyc:G_\Q\to\Z_p^\times$ is $p$-adic cyclotomic character. We then get automatically a universal weight character (cf. remark \ref{lambdaadicHeckeop} for the notation):
\[
\langle\cyc\rangle_\Lambda:G_\Q\to\Lambda^\times\qquad g\mapsto [\langle\cyc(g)\rangle]=\langle\,\cdot\,\rangle_\Lambda\circ\cyc(g)
\]
and, for $(R,\varphi)\in\CCC_\Lambda$, we set $\langle\cyc\rangle_R=\varphi\circ\langle\cyc\rangle_\Lambda =\langle\,\cdot\,\rangle_R\circ\cyc:G_\Q\to R^\times$.

Since we assume that $p$ does not divide the class number of $K$, we have $W_K=1+p\Z_{p^2}$ (cf. lemma \ref{Wlambda}) and moreover (with the notation of remark \ref{Zp2remark})
\begin{equation}
\label{comparecyclambda}
\langle\cyc\rangle|_{G_K}=\langle\lambda\rangle\cdot\langle\lambda\rangle^\sigma.
\end{equation}
\end{rmk}

By the work of Hida and Wiles, it is known that one can attach to $\Hf$ a \emph{big} Galois representation $\V_{\!\Hf}$, which can be realized as a free module of rank $2$ over $\Lambda_\Hf[1/p]$ equipped with a continuous action of $G_\Q$, specializing for all $x\in\W_{\Lambda_\Hf}^{cl}$ to the dual $V_p(\Hf_{\!x})$ of the $p$-adic Deligne representation attached to $\Hf_{\!x}$ (or, in case $\Hf_{\!x}$ is the $p$-stabilization of a newform of level $N_\Hf$, to the dual of the representation attached to such newform). In particular it holds that $\det(\V_{\!\Hf})=\omega_\mathrm{cyc}\cdot\langle\cyc\rangle_{\Lambda_\Hf}$. We refer to \cite[section 5]{BSVast1} for a detailed discussion concerning such Galois modules.

\medskip
We defined a big Galois representation $\V_{\!\Hg}$ (resp. $\V_{\!\Hh}$) attached to $\Hg$ (resp. $\Hh$) as
\[
\V_{\!\Hg}=\Ind_K^\Q\,\pmb{\eta}_1\qquad \text{(resp.}\V_{\!\Hh}=\Ind_K^\Q\,\pmb{\eta}_2\text{),}
\]
where $\pmb{\eta}_1$ (resp. $\pmb{\eta}_2$) is the big Galois character valued in $\Lambda_\mathrm{Hida}[1/p]$ constructed as in section \ref{alahida}.

\begin{nota}
We will write $R_K:=\Lambda_\mathrm{Hida}$ in what follows, to simplify the notation. We will also write $\langle\lambda\rangle_{R_K}:G_K\to R_K^\times$ for the big Galois character given by $g\mapsto [\langle\lambda(g)\rangle]$.
\end{nota}

\begin{lemma}
We have
\[
\det(\V_{\!\Hg})=\varepsilon_K\cdot\eta_1^\mathrm{cen}\cdot\langle\cyc\rangle\cdot\langle\cyc\rangle_{R_K},\qquad \det(\V_{\!\Hh})=\varepsilon_K\cdot\eta_2^\mathrm{cen}\cdot\langle\cyc\rangle\cdot\langle\cyc\rangle_{R_K}.
\]
\begin{proof}
It follows easily from equation \eqref{comparecyclambda}.
\end{proof}   
\end{lemma}

\medskip
Consider the Galois representation $\V:=\V_{\!\Hf}\hat{\otimes}_L\V_{\!\Hg}\hat{\otimes}_L\V_{\!\Hh}$. It is a free $\RR:=R_{\Hf\Hg\Hh}[1/p]$-module of rank $2$ and it follows immediately from the above discussion and our assumptions that
\[
\det(\V)=\omega_\mathrm{cyc}^6\cdot\cyc^{-3}\cdot(\langle\cyc\rangle_{\Lambda_\Hf}\hat{\otimes}\langle\cyc\rangle_{R_K}\hat{\otimes}\langle\cyc\rangle_{R_K}).
\]
Since $p$ is odd, there exists a character $\pmb{\chi}_{\Hf\Hg\Hh}=\pmb{\chi}:G_\Q\to\RR^\times$ such that $\cyc\cdot \pmb{\chi}^2=\det(\V)$, i.e., we can write
\[
\pmb{\chi}=\omega_\mathrm{cyc}\cdot\langle\cyc\rangle^{-2}\cdot(\langle\cyc\rangle_{\Lambda_\Hf}\hat{\otimes}\langle\cyc\rangle_{R_K}\hat{\otimes}\langle\cyc\rangle_{R_K})^{1/2}.
\]
If we define
\[
\V^\dagger:=\V\otimes_\RR\RR(\pmb{\chi}^{-1})\,,
\]
then one checks easily that such representation is Kummer self-dual, i.e.
\[
(\V^\dagger)^\vee(1)=\Hom_\RR(\V^\dagger,\RR)(1)\cong\V^\dagger.
\]
We want to study the specializations $\V^\dagger(w)$ for a suitable $w\in\Omega_{\Hf\Hg\Hh}$.

\begin{defi}
We define the following big Galois characters
\begin{equation}
\label{bigphipsi}    \pmb{\varphi}:=\eta_1\eta_2\langle\lambda\rangle^\sigma\langle\lambda\rangle^{-1}\cdot\pmb{\lambda}_\mathrm{ac}\qquad\text{and}\qquad \pmb{\psi}:=\eta_1\eta_2^\sigma\cdot\pmb{\lambda}_\mathrm{ac}\,,
\end{equation}
where
\[
\pmb{\lambda}_\mathrm{ac}:G_K\to R_K^\times\qquad \pmb{\lambda}_\mathrm{ac}:=\langle\lambda\rangle_{R_K}^{1/2}\cdot (\langle\lambda\rangle^\sigma_{R_K})^{-1/2}\,.
\]
\end{defi}

\begin{rmk}
Note that $\W_{R_K}(\C_p)\cong\Hom^{cont}_{grp}(W_K,\C_p^\times)$ has a natural group structure, so it makes sense to multiply or invert weights.    
\end{rmk}

\begin{lemma}
\label{Galoisfact}
Let $w=(x,y,z)\in\Omega_{\Hf\Hg\Hh}$ with $k=k_x$ even and let $f^\circ$ be the newform associated to $\Hf_{\!x}$ (as in proposition \ref{classinterp}). Then there is a decomposition
\[
\V^\dagger(w)\cong\left(\left(V_p(\Tilde{f^\circ})\otimes_{L[w]}\Ind_K^\Q\pmb{\varphi}_{y\cdot z}\right)\oplus\left(V_p(\Tilde{f^\circ})\otimes_{L[w]}\Ind_K^\Q\pmb{\psi}_{y/z}\right)\right)\left(-\tfrac{k}{2}\right)\,,
\]
where $\Tilde{f^\circ}:=f^\circ\otimes\omega^{k/2-1}\varepsilon_x^{-1/2}$.

Moreover, setting $l=k_y$ and $m=k_z$, the Hecke character of $K$ attached to $\pmb{\varphi}_{y\cdot z}$ (resp. to $\pmb{\psi}_{y/z}$) is anticyclotomic and has $\infty$-type $(\frac{l+m-2}{2},\frac{2-l-m}{2})$ (resp. $(\frac{l-m}{2},\frac{m-l}{2})$).
\begin{proof}
This is an easy computation. The only passage when one has to be slightly careful consists in observing that for a $G_K$-character $\eta$ and an even $G_\Q$-character $\chi$ it holds $\left((\Ind_K^\Q(\eta)\right)(\chi)=\Ind_K^\Q(\eta\cdot\chi|_{G_K})$.
\end{proof}
\end{lemma}

\subsection{Improvement of the triple product \texorpdfstring{$p$}{p}-adic \texorpdfstring{$L$}{L}-function in our setting}
\label{improvedLfctsection}
We let $M:=c^2\cdot d_K\cdot N_\Hf$ in what follows.

\medskip
Inspired by the level adjustment performed by Hsieh in \cite[section 3.4]{Hsi2021}, we are led to consider the following \emph{test vectors} associated to our families $\Hg$ and $\Hh$, namely we set
\begin{equation}
\label{leveladj}
    \Hg^*:=\Hg(q^{N_\Hf})\in\SSS_{\Omega_\mathrm{Hida}}(M,\chi_\Hg,R_K)\,,\qquad \Hh^*:=\Hh(q^{N_\Hf})\in\SSS_{\Omega_\mathrm{Hida}}(M,\chi_\Hh,R_K).
\end{equation}
One can check that our adjustment matches Hsieh's more general version, in view of the following facts concerning the local automorphic types for the specializations of the families $\Hf$, $\Hg$ and $\Hh$.
\begin{prop}
Let $\ell$ be a prime different from $p$. Let $w=(x,y,z)\in\Omega_{\Hf\Hg\Hh}$ and write $(f,g,h)=(\Hf_{\!x},\Hg_y,\Hh_z)$. Denote by $\pi_{\xi,\ell}$ the local component at $\ell$ of the automorphic representation $\pi_{\xi}$ attached to $\xi\in\{f,g,h\}$. Then the following facts hold.
\begin{itemize}
    \item [(i)] The automorphic type of $\pi_{\xi,\ell}$ does not depend on the chosen specialization for $\xi\in\{f,g,h\}$ (rigidity of automorphic types).
    \item [(ii)] If $\ell\nmid M$, then $\pi_{\xi,\ell}$ is an unramified principal series representation for $\xi\in\{f,g,h\}$.
    \item [(iii)] If $\ell\mid N_\Hf$, then $\pi_{f,\ell}$ is special, while $\pi_{g,\ell}$ and $\pi_{h,\ell}$ are unramified principal series.
    \item [(iv)] If $\ell\mid c^2d_K$, then $\pi_{f,\ell}$ is an unramified principal series, while $\pi_{g,\ell}$ and $\pi_{h,\ell}$ are ramified principal series.
\end{itemize}
\begin{proof}
All the assertions regarding $\pi_{f,\ell}$ are well-known for Hida families and for the choice of squarefree tame level $N_\Hf$ and trivial character in our setting. The assertions regarding $\pi_{g,\ell}$ and $\pi_{h,\ell}$ follow from the explicit description of the Weil-Deligne representations which correspond to them via the local Langlands correspondence. Here we use the assumption that $d_K$ is odd and that the prime divisors of $c$ split in $K$ to grant that the restriction of $V_p(g)$ and $V_p(h)$ to a decomposition group at $\ell$ is reducible when $\ell\mid d_K$.
\end{proof}
\end{prop}

Along the lines of \cite[proposition 6.12]{Hsi2021}, we can thus define the so-called \emph{fudge factors} at the primes dividing $M$.

\begin{prop}
For each $\ell\mid M$, there exists a unique element $\fid_{\Hf\Hg\Hh, \ell}\in R_{\Hf\Hg\Hh}^\times$ such that for all $w\in\Omega_{\Hf\Hg\Hh}$ it holds
\[
(\fid_{\Hf\Hg\Hh, \ell})_w=\mathscr{I}^*_{\Pi_w,\ell}\,,
\]
with $\mathscr{I}^*_{\Pi_w,\ell}$ as in proposition \ref{classinterp}.
\begin{proof}
This is proven (adapting Hsieh's methods) in the same way as in \cite[section 5.1]{Fuk2022}.
\end{proof}
\end{prop}

\begin{defi}
\label{improvedLfct}
We define the element
\[
\mathcal{L}_p^f(\Hf,\Hg,\Hh):=\Ls^f_p(\Hf,\Hg^*,\Hh^*)\cdot\prod_{\ell\mid M}\fid_{\Hf\Hg\Hh, \ell}^{-1/2}\in R_{\Hf\Hg\Hh}
\]
and call it the square root $f$-unbalanced $p$-adic triple product $L$-function attached to our triple $(\Hf,\Hg,\Hh)$.   
\end{defi}

\begin{cor}
With the above notation, for all $w\in\Omega_{\Hf\Hg\Hh}$ lying in the $\Hf$-unbalanced region, it holds
\begin{equation}
\label{interpupgrade}
    (\mathcal{L}_p^f(\Hf,\Hg,\Hh)(w))^2= \frac{L^*(\Pi,1/2)}{\zeta_\Q(2)^2\cdot\Omega_f^2}\cdot\mathscr{I}^{unb}_{\Pi_w,p}\,.
\end{equation}

\begin{proof}
Obvious from the formula \eqref{classinterp} and the definition of $\mathcal{L}_p^f(\Hf,\Hg,\Hh)$.
\end{proof}
\end{cor}

\bigskip
We are left to find a more explicit description of the local integral $\mathscr{I}^{unb}_{\Pi_w,p}$. We will fix a triple of weights $w=(x,y,z)\in\Omega_{\Hf\Hg\Hh}$ which is $\Hf$-unbalanced. Write $k=k_x,l=k_y,m=k_z$ as usual, so that $k\geq l+m$. Assume furthermore that $k$ is even.

Let $(f,g,h)=(\Hf_{\!x},\Hg_y,\Hh_z)$ as above and, only for this section, set
\[
\pi_1:=\pi_{f,p}\otimes\tilde{\chi}_1\, ,\qquad\pi_2:=\pi_{g,p}\, ,\qquad\pi_3:=\pi_{h,p}.
\]
where
\[
\tilde{\chi}_1=\omega^{(k+l+m-6)/2}\cdot(\varepsilon_x\varepsilon_y\varepsilon_z)^{-1/2}
\]
Let $\chi_1=\alpha_{f,p}\cdot\tilde{\chi}_1$, where $\alpha_{f,p}$ denotes the unramified character of $\Q_p^\times$ such that $\alpha_{f,p}(p)=a_p(f)p^{(1-k)/2}$.

\medskip
Then $\pi_i$ is an irreducible smooth representations of $\GL_2(\Q_p)$ for $i=1,2,3$ and, by our assumptions, we know the following.

\begin{lemma}
\label{localpauttypes}
The representations $\pi_2$ and $\pi_3$ are always supercuspidal. The representation $\pi_1$ satisfies one of the following:
\begin{itemize}
    \item [(a)] $\pi_1$ is the principal series $\pi_1=\chi_1\boxplus\nu_1$ with $\nu_1=\omega_{f,p}\chi_1^{-1}$ where $\omega_{f,p}$ is the $p$-component of the central character of $\pi_{f,p}$;
    \item [(b)] $\pi_1$ is the special representation $\pi_1=\chi_1|\cdot|^{-1/2}\mathrm{St}$.
\end{itemize}
The latter case happens if and only if $x=2$ and $f=\Hf_{\!2}$ is $p$-new.
\begin{proof}
All the assertions concerning $\pi_1$ are well-known for Hida families. The fact that $\pi_2$ and $\pi_3$ are always supercuspidal follows from the fact that $g$ and $h$ are theta series attached to a Hecke character of $K$ ramified at $p$ (recall that the prime $p$ is inert in $K$ by assumption).
\end{proof}
\end{lemma}

\begin{prop}
In the above setting, we have that
\[
\mathscr{I}^{unb}_{\Pi_w,p}=\frac{L(\pi_2\otimes\pi_3\otimes\chi_1, 1/2)}{\varepsilon(\pi_2\otimes\pi_3\otimes\chi_1,1/2)\cdot L(\pi_2\otimes\pi_3\otimes\nu_1,1/2)\cdot L(\pi_1\otimes\pi_2\otimes\pi_3,1/2)}\,.
\]
\begin{proof}
This follows adapting \cite[proposition 5.4]{Hsi2021} in the same way as it is suggested in \cite[remark 3.4.7]{Fuk2022}. With the notation of \cite[section 5.3.2]{Hsi2021}, this means that for the calculation of $\mathscr{I}^{unb}_{\Pi_w,p}$ we choose as test vector
\[
\phi_p^*=(f_1^\ord)^0\otimes W_2\otimes\theta_p^\mathbbm{k}W_3\,,
\]
where, for $i=2,3$, $W_i$ is the local Whittaker newform for $\pi_i$, as defined in \cite[section 2.4.5]{Hsi2021} (instead of the normalized ordinary Whittaker function chosen for the computation in \cite[proposition 5.4]{Hsi2021}). Since for $i=2,3$ it still holds that
\[
W_i\left(\begin{pmatrix}
    a & 0\\
    0 & 1
\end{pmatrix}\right) = \mathbbm{1}_{\Z_p^\times}(a)\,\quad\text{for every }a\in\Q_p^\times\,,
\]
the same calculation of \cite[page 488]{Hsi2021} for the local Rankin-Selberg integral (denoted $\Psi(W_2,\theta_p^\mathbbm{k}W_3,\rho(t_n)f_1^\ord)$ in \cite{Hsi2021}) yields the desired formula.
\end{proof}
\end{prop}

We can give an even more explicit description of $\mathscr{I}^{unb}_{\Pi_w,p}$. Write $\varphi$ (resp. $\psi$) to denote - again only in this section - the $p$-component of $\pmb{\varphi}_{y\cdot z}$ (resp. $\pmb{\psi}_{y/z}$) seen as Hecke character of $K$. Let also $\mu_1$ and $\mu_2$ denote the characters of $\Q_{p^2}$ given by
\[
\mu_1=(\alpha_{f,p}\cdot\omega_{f,p}^{-1/2})\circ N_{\Q_{p^2}/\Q_p}\, ,\qquad
\mu_2=(\alpha_{f,p}^{-1}\cdot\omega_{f,p}^{-1/2})\circ N_{\Q_{p^2}/\Q_p}\,
\]
and set $\pi_1'=\pi_{f,p}\otimes\omega_{f,p}^{-1/2}$.

\begin{prop}
\label{localfactp}
With the above notation, it holds $\mathscr{I}^{unb}_{\Pi_w,p}=\mathscr{I}_{\varphi,w}\cdot\mathscr{I}_{\psi,w}$, where for $\eta\in\{\varphi,\psi\}$ we set
\begin{equation}
    \mathscr{I}_{\eta,w}:=\frac{L(\pi(\eta\mu_1), 1/2)}{\varepsilon(\pi(\eta\mu_1),1/2)\cdot L(\pi(\eta\mu_2),1/2)\cdot L(\pi_1'\otimes\pi(\eta),1/2)}\, .
\end{equation}
Moreover one can compute $\mathscr{I}_{\eta,w}$ as follows.
\begin{itemize}
    \item [(1)] Assume that we are in case (a) of lemma \ref{localpauttypes} and that the character $\eta\mu_1$ is unramified, then
    \[
    \mathscr{I}_{\eta,w}=\left(1-\frac{p^{k-2}}{a_p(f)^2}\right)^2\, .
    \]
    \item [(2)] Assume that we are in case (b) of lemma \ref{localpauttypes} and that the character $\eta\mu_1$ of $\Q_{p^2}^\times$ is unramified, then
    \[
    \mathscr{I}_{\eta,w}=1-\frac{p^{k-2}}{a_p(f)^2}=1-a_p(f)^{-2}\, .
    \]
    \item [(3)] Assume that the character $\eta\mu_1$ of $\Q_{p^2}^\times$ is ramified of level $n$, then
    \[
    \mathscr{I}_{\eta,w}=\left(\frac{p}{a_p(f)^2}\right)^n\cdot\frac{p^{n(k-2)}}{W(\tilde{\eta})}\, ,
    \]
    where $\tilde{\eta}$ is the unitary character of $\Q_{p^2}$ given by $\eta\mu_1$ on $\Z_{p^2}^\times$ and such that $\tilde{\eta}(p)=1$ and $W(\tilde{\eta})$ denotes the root number of $\tilde{\eta}$, defined as
    \[
    W(\tilde{\eta})=\varepsilon(\tilde{\eta},1/2)\, ,
    \]
    which is an algebraic integer of complex absolute value $1$.
\end{itemize}
\begin{proof}
The factorization $\mathscr{I}^{unb}_{\Pi_w,p}=\mathscr{I}_{\varphi,w}\cdot\mathscr{I}_{\psi,w}$ follows directly from the corresponding factorization at the level of Galois representations given in lemma \ref{Galoisfact} and the local Langlands correspondence for $\GL_2(\Q_p)$. Hence we know that
\[
\mathscr{I}_{\eta,w}=\frac{L(\eta\mu_1, 1/2)}{\varepsilon(\eta\mu_1,1/2)\cdot L(\eta\mu_2,1/2)\cdot L(\pi_1'\otimes\pi(\eta),1/2)}\, .
\]
Now note that for $\eta\in\{\varphi,\psi\}$ we have that $\eta\mu_1$ is a unitary character $\Q_{p^2}^\times\to\C^\times$, since $\pmb{\varphi}_{y\cdot z}$ and $\pmb{\psi}_{y/z}$ are anticyclotomic and $\mu_i$ is unitary for $i=1,2$. The fact that $\pmb{\varphi}_{y\cdot z}$ and $\pmb{\psi}_{y/z}$ are anticyclotomic also implies that $\eta(p)=1$.

\medskip
We can proceed depending on the three cases, applying the known facts from Tate's thesis for the definition local $L$-factors and $\varepsilon$-factors attached to Hecke characters.

\medskip
\begin{itemize}
    \item [(1)] If $\eta\mu_1$ is unramified, then $\varepsilon(\eta\mu_1,1/2)=1$. Moreover, if $\pi_1'$ is the unramified principal series $\pi_1'=\alpha_{f,p}\cdot\omega_{f,p}^{-1/2}\boxplus\alpha_{f,p}^{-1}\cdot\omega_{f,p}^{-1/2}$, then $L(\pi_1'\times\pi(\eta),s)=L(\eta\mu_1,s)\cdot L(\eta\mu_2,s)$. Hence
    \[
    \mathscr{I}_{\eta,w}=L(\eta\mu_2,1/2)^{-2}=\left(1-\eta\mu_2(p)\cdotp^{-1}\right)^2=\left(1-\frac{p^{k-2}}{a_p(f)^2}\right)^2\,.
    \]
    \item [(2)] If $\pi_1'=\alpha_{f,p}\cdot\omega_{f,p}^{-1/2}|\cdot|^{-1/2}\mathrm{St}$, then $L(\pi_1'\times\pi(\eta),s)=L(\eta\mu_1,s)$. Hence
    \[
    \mathscr{I}_{\eta,w}=L(\eta\mu_2,1/2)^{-1}=1-\eta\mu_2(p)p^{-1}=1-a_p(f)^{-2}\,,
    \]
    where we used that this situation can only occur with $x=k=2$.
    \item [(3)] If $\eta\mu_1$ is ramified of level $n$ (so that necessarily also $\eta\mu_2$ is ramified), all the $L$-factors involved are equal to $1$, so that $\mathscr{I}_{\eta,w}=\varepsilon(\eta\mu_1,1/2)^{-1}$ and by Tate's thesis we know $ \varepsilon(\eta\mu_1,1/2)=\eta\mu_1(p)^n\cdot\varepsilon(\tilde{\eta},1/2)$. Hence
    \[
    \mathscr{I}_{\eta,w}=\varepsilon(\eta\mu_1,1/2)^{-1}=\eta\mu_1(p)^{-n}\cdot W(\tilde{\eta})^{-1}=\left(\frac{p}{a_p(f)^2}\right)^n\cdot\frac{p^{n(k-2)}}{W(\tilde{\eta})}\,.
    \]
\end{itemize}
\end{proof}
\end{prop}

\begin{rmk}
We observe that the results of the above computation match perfectly the shape of the modification of the Euler factor at $p$ (for the Galois theoretic side) described in \cite[pagg. 162-163]{Coa1991}, also in the cases of \emph{bad reduction} at $p$.
\end{rmk}

We have some control on the root numbers appearing in proposition \ref{localfactp} (case (3)).
\begin{lemma}
\label{rootnumberlemma}
With the notation introduced above, if $x=k\equiv 2\mod (p-1)$ and the character $\eta\in\{\varphi,\psi\}$ is ramified, then $W(\tilde{\eta})=W(\eta)\in\{\pm 1\}$. Moreover the sign $W(\varphi)$ (resp. $W(\psi)$) depends only on the parity of $j_1=(l+m-2)/2$ (resp. $j_2=(l-m)/2$).
\begin{proof}
Note that under our assumptions the character denoted $\mu_1$ above is unramified and $\eta=\tilde{\eta}$ is of finite order and trivial on $\Q_p^\times$. We can thus apply \cite[proposition 3.7]{MS2000} to a suitable twist of $\eta$ to deduce that $W(\eta)=\eta^{-1}(\alpha)$, where $\alpha\in \Q_{p^2}^\times$ is a primitive $2(p-1)$-th root of unit, so that $1=\eta(-1)=\eta(\alpha)^{-2}$. In particular this shows that $W(\eta)\in\{\pm 1\}$.

Recall that $\Z_{p^2}^\times=\mu_{p^2-1}\times (1+p\Z_{p^2})$. Thus the only way one can affect the sign $W(\eta)$ is changing the weights $l,m$. More precisely, one can check (cf. remark \ref{cyclolambda}) that
\[
\varphi|_{\mu_{p^2-1}}=\eta_1\eta_2|_{\mu_{p^2-1}}\cdot (-)^{\tfrac{(p-1)(l+m-2)}{2}}\, ,\qquad \psi|_{\mu_{p^2-1}}=\eta_1\eta_2^\sigma|_{\mu_{p^2-1}}\cdot (-)^{\tfrac{(p-1)(l-m)}{2}}\, .
\]
Writing $\alpha=\zeta^{(p+1)/2}$ for $\zeta$ a primitive $(p^2-1)$-th root of $1$, we see that the sign $W(\varphi)$ (resp. $W(\psi)$) depends only on the parity of $j_1=(l+m-2)/2$ (resp. $j_2=(l-m)/2$).
\end{proof}
\end{lemma}
\subsection{Anticyclotomic \texorpdfstring{$p$}{p}-adic \texorpdfstring{$L$}{L}-functions}
\medskip
As in the introduction, let $H_n$ denote the ring class field of $K$ of conductor $cp^n$ and let $H_\infty$ be the union of all the $H_n$'s. It follows that the big characters $\pmb{\varphi}$ and $\pmb{\psi}$ (defined in equation \ref{bigphipsi}) factor through $\mathscr{G}_\infty:=\Gal(H_\infty/K)$. With the same notation as in remark \ref{Zp2remark}, we can identify $\Gamma^-=\Gal(K_\infty^-/K)$ (the Galois group of the anticyclotomic $\Z_p$-extension of $K$) with the maximal $\Z_p$-free quotient of $\mathscr{G}_\infty$, i.e., there is an exact sequence
\[
0\to\Delta_c\to\mathscr{G}_\infty\to\Gamma^-\to 0
\]
of abelian groups with $\Delta_c$ a finite group and $\Gamma^-\cong\Z_p$. We fix a non-canonical isomorphism $\mathscr{G}_\infty\cong\Delta_c\times\Gamma^-$ once and for all. Notice that $\pmb{\lambda}_\mathrm{ac}$ will factor through $\Gamma^-$.

\medskip
As in lemma \ref{rootnumberlemma}, set $j_1:=\tfrac{l+m-2}{2}$ and $j_2:=\tfrac{l-m}{2}$. If we assume moreover that the triple of weights $w=(k,y,z)$ is $\Hf$-unbalanced (i.e. $k\geq l+m$), then $|j_i|<\frac{k}{2}$ for $i=1,2$.

\medskip
Building up on previous work of Bertolini-Darmon (\cite{BD1996}, \cite{BD1998}) and Chida-Hsieh (\cite{CH2018}), Castella and Longo in \cite{CL2016} have constructed so-called \emph{big theta elements}, denoted 
\begin{equation}
    \Theta_\infty^{\mathrm{Heeg}}(\Hf)\in R_{\Hf,\Gamma^-}:=\Lambda_\Hf\hat{\otimes}_{\OO_L}\OO_L[\![\Gamma^-]\!]
\end{equation}
attached to the Hida family $\Hf$ and the quadratic imaginary field $K$ (satisfying a suitable Heegner hypothesis relative to the tame level of $\Hf$). The two variables are given by the weight specializations for $\Hf$ and by continuous characters $\hat{\nu}:\Gamma^-\to\C_p^{\times}$ such that the associated algebraic Hecke character $\nu:\Ad_K^\times/K^\times\to\C^\times$ has infinity type $(j,-j)$ with $|j|<k/2$. We let $\mathfrak{X}^{\mathrm{crit}}_{p,k}$ to denote the set of characters $\hat{\nu}$ satisfying such requirement for a fixed $k$. The specializations of the square of $\Theta_\infty^{\mathrm{Heeg}}(\Hf)$ at $(k,\hat{\nu})$ with $k\geq 2$ even integer and $\hat{\nu}\in\mathfrak{X}^{\mathrm{crit}}_{p,k}$ interpolate the (algebraic part of the) special values $L(\Hf_{\!k}^\circ/K,\nu,k/2)$.

\medskip
Following the strategy of Castella and Longo applied to the more general construction of Hung (\cite{Hu2017}), one can construct a big theta element $\Theta_\infty^{\mathrm{Heeg}}(\Hf,\chi_t)\in R_{\Hf,\Gamma^-}$ associated with the Hida family $\Hf$ and a branch character $\chi_t$ of conductor $c$ (i.e. a character of the finite group $\Delta_c$).

\begin{rmk}
\label{anticycchoicesrmk}
The construction of $\Theta_\infty^{\mathrm{Heeg}}(\Hf,\chi_t)$ depends on the following choices that we fix from now on:
\begin{itemize}
    \item[(a)] a factorization $N_\Hf^+\OO_K=\mathfrak{N}^+\cdot\overline{\mathfrak{N}^+}$,
where recall that $N_\Hf^+$ is the product of the prime divisors of $N_\Hf$ that split in $K$;
\item[(b)] a family of quaternionic modular forms $\boldsymbol{\Phi}$ associated to $\Hf$, with the property that there exists an open neighbourhood $U_\Hf$ of $2$ in $\W_{\Lambda_\Hf}(\OO_L)$ such that for all $k\in U_{\Hf}\cap \Z_{\geq 2}$ it holds
\[
\boldsymbol{\Phi}_k=\lambda_{B,k}\cdot\varphi_k\,,
\]
where $\lambda_{B,k}\in L^\times$ and $\varphi_k$ corresponds to $\Hf_{\!k}^\circ$ via a version of the Jacquet-Langlands correspondence. 
\end{itemize}

We can (and will) choose the following normalizations for $\boldsymbol{\Phi}$:
\begin{itemize}
    \item [(i)] $\lambda_{B,2}=1$;
    \item [(ii)] $\eta_{\Hf_{\!k}^\circ,N^-}=1$ for $k\in U_\Hf\cap\Z_>2$.
\end{itemize}
The period $\eta_{\Hf_{\!k}^\circ,N^-}$ (appearing in the following proposition) is defined as a suitable Petersson norm of $\varphi_k$, which we can normalize to be $1$ (this will determine $\varphi_k$ up to sign). We refer to \cite[theorem 2.5]{BD2007} for the existence of $\boldsymbol{\Phi}$ and its properties and to \cite[equations 3.9 and 4.3]{CH2018} for the description of $\eta_{\Hf_{\!k}^\circ,N^-}$ as Petersson norm (Chida-Hsieh's notation is $\langle f_{\pi'}, f_{\pi'}\rangle_R$).
\end{rmk}

\begin{prop}
\label{anticyclotomicinterpolationprop}
Fix an even integer $k\in U_\Hf\cap \Z_{\geq 2}$ and a character $\hat{\nu}\in\mathfrak{X}^{\mathrm{crit}}_{p,k}$ of conductor $p^n$. Write $f=\Hf_{\!k}$ and $f^\circ=\Hf_{\!k}^\circ$ (with the usual conventions). Then:
\begin{equation}
\label{anticycinterp}
\left(\Theta_\infty^{\mathrm{Heeg}}(\Hf/K,\chi_t)\right)^2(k,\hat{\nu})=\lambda_B(k)^2\cdot C_p(f,\chi_t\nu)\cdot e_p(f,\chi_t\nu)\cdot\frac{L(f^\circ/K,\chi_t\nu,k/2)}{\Omega_{f^\circ,N^-}}
\end{equation}
where:
\begin{itemize}
    \item [(i)] setting $u_K=\tfrac{\#\OO_K^\times}{2}$ and $\delta_K:=\sqrt{d_K}$, one has
    \[
    C_p(f,\chi_t\nu):=(-1)^{\frac{2+2j-k}{2}}\cdot\Gamma(k/2+j)\cdot\Gamma(k/2-j)\cdot c\cdot \delta_K^{k-1}\cdot u_K^2\cdot\varepsilon(\pi_{f,p},1/2)\cdot\chi_t\nu(\mathfrak{N}^+) ;
    \]
    \item [(ii)] 
    \[
    e_p(f,\chi_t\nu)=\left\{
    \begin{aligned}
    &\left(\frac{p}{a_p(f)^2}\right)^n\cdot p^{n(k-2)} && \text{ if }n>0\\
    & \left(1-\frac{p^{k-2}}{a_p(f)^2}\right)^2 && \text{ if $n=0$ and $f$ is $p$-old}\\
    & 1-\frac{p^{k-2}}{a_p(f)^2} && \text{ if $n=0$ and $f$ is $p$-new}
    \end{aligned}
    \right. \quad ;
    \]
    \item [(iii)] $\Omega_{f^\circ,N^-}$ is Gross's period, that we can write as
    \begin{equation}
        \Omega_{f^\circ,N^-}=\frac{(4\pi)^k\cdot\|f^\circ\|_{Pet}^2\cdot\zeta_\Q(2)\cdot[\SL_2(\Z):\Gamma_0(N_\Hf)]}{\eta_{f^\circ,N^-}}\; .
    \end{equation}

\end{itemize}
\begin{proof}
This follows from the work of Chida-Hsieh \cite{CH2018}, Hung \cite{Hu2017} and Castella-Longo \cite{CL2016}. We refer to \cite[section 4.2]{CL2016} and to \cite[theorem 5.6]{Hu2017} for the interpolation formula. Gross's period is defined as in \cite[equation (5.2)]{Hu2017} (translated into our notation).
\end{proof}
\end{prop}

\begin{rmk}
\label{definitenonvanish}
We keep the notation of proposition \ref{anticyclotomicinterpolationprop}. The Heegner hypothesis (iii) on $N_\Hf$ in assumption \ref{tripleprodass} implies that the sign of the functional equation for $L(f^\circ/K,\chi,s)$ is $+1$ for every anticyclotomic Hecke character $\chi$ of $K$ of conductor coprime to $N_\Hf\cdot d_K$ (unless $k=2$, $\Hf_{\!2}$ is $p$-new and $\chi$ is unramified at $p$), i.e. we are in the so-called \emph{definite setting}. One of the main results of \cite{Hu2017} (namely theorem C in the introduction), generalizing work of Vatsal \cite{Va2022} and Chida-Hsieh \cite{CH2018}, implies that in our setting it holds $L(f^\circ/K,\chi_t\nu,k/2)\neq 0$ for all but finitely many $\hat{\nu}\in\mathfrak{X}^\mathrm{crit}_{p,k}$.
\end{rmk}

\subsection{Factorization of the triple product \texorpdfstring{$p$}{p}-adic \texorpdfstring{$L$}{L}-function}
We consider the automorphism $s$ of $R_K\hat{\otimes}_{\OO_L}R_K$ in $\CCC_{\OO_L}$ given by the assignment
\[
[\gamma]\otimes[\delta]\mapsto[\gamma^{1/2}\delta^{1/2}]\otimes[\gamma^{1/2}\delta^{-1/2}]
\]
on group-like elements (note that again it is important that $p\neq 2$ for this to be a well-defined automorphism). 

Let again $K_\infty$ denote the (unique) $\Z_p^2$-extension of $K$. Recall (remark \ref{Zp2remark}) that the character $\langle\lambda\rangle$ induces an isomorphism $\Gamma_\infty\cong W_K$. The natural projection $\Gamma_\infty\twoheadrightarrow\Gamma^-$ can be described as $\gamma\mapsto\gamma^{1/2}(\gamma^\sigma)^{-1/2}$. Accordingly, we get a morphism
\begin{equation}
\label{anticycproj}
    \tau:R_K\twoheadrightarrow\OO_L[\![\Gamma^-]\!]\,.
\end{equation}

\begin{nota}
We set $\varphi_t:=\eta_1\eta_2|_{\Delta_c}$ and $\psi_t:=\eta_1\eta_2^\sigma|_{\Delta_c}$. With respect to the chosen isomorphism $\mathscr{G}_\infty\cong\Delta_c\times\Gamma^-$, we also define the characters of $\Gamma^-$ given by $\varphi^-:=\eta_1\eta_2|_{\Gamma^-}$ and $\psi^-:=\eta_1\eta_2^\sigma|_{\Gamma^-}$.   
\end{nota}

Note that the assigments $[\gamma]\mapsto\varphi^-(\gamma)[\gamma]$ (resp. $[\gamma]\mapsto\psi^-(\gamma)[\gamma]$) define $\OO_L$-linear automorphisms $\varphi^-:\OO_L[\![\Gamma^-]\!]\cong\OO_L[\![\Gamma^-]\!]$ (resp. $\psi^-:\OO_L[\![\Gamma^-]\!]\cong\OO_L[\![\Gamma^-]\!]$), since $|\varphi^-(\gamma)-1|_p<1$ (resp. $|\psi^-(\gamma)-1|_p<1$) for $\gamma\in\Gamma^-$. By slight abuse of  denote by $\varphi^-$ (resp. $\psi^-$) the automorphism of $R_{\Hf,\Gamma^-}$ given by the identity on $\Lambda$ and $\varphi^-$ (resp. $\psi^-$) on $\OO_L[\![\Gamma^-]\!]$.

\begin{lemma}
\label{weighttransitionlemma}
Consider the composition
\[
\mathrm{pr}_{ac}:R_{\Hf\Hg\Hh}\xrightarrow[\cong]{1\otimes s}R_{\Hf\Hg\Hh}\xrightarrow{1\otimes\tau\otimes\tau}\Lambda_\Hf\hat{\otimes}_{\OO_L}\OO_L[\![\Gamma^-]\!]\hat{\otimes}_{\OO_L}\OO_L[\![\Gamma^-]\!]\,.
\]
Given a specialization $(k,\hat{\nu},\hat{\mu})\in\W^{cl}_{\Lambda_\Hf,\Z}\times\mathfrak{X}^\mathrm{crit}_{p,k}\times\mathfrak{X}^\mathrm{crit}_{p,k}$ (with $k\geq 2$ even integer), then the specializations in $\Omega_{\Hf\Hg\Hh}$ which lift $(k,\hat{\nu},\hat{\mu})$ are $\Hf$-unbalanced triples $w=(k,y,z)$ with the property that
\begin{equation}
\label{condyz}
    \hat{\nu}=(yz)|_{\Gamma^-}\cdot\langle\lambda\rangle^\sigma\langle\lambda\rangle^{-1}\; ,\qquad \hat{\mu}=(y/z)|_{\Gamma^-}\; .
\end{equation}
Moreover, we can always find such $y\in\Omega_\Hg$ and $z\in\Omega_\Hh$ for given $\hat{\nu}$ and $\hat{\mu}$ such that $w=(k,y,z)$ is $\Hf$-unbalanced.
\begin{proof}
This is an easy exercise.
\end{proof}
\end{lemma}

\begin{nota}
\label{R-notation}
Now let $\sigma_{\mathfrak{N}^+}\in\mathscr{G}_\infty$ denote the projection to $\mathscr{G}_\infty$ of the element of $G_K$ corresponding to $\mathfrak{N}^+$ by class field theory. We write $(\sigma_c,\gamma_{\mathfrak{N}^+}^{-2}):=\sigma_{\mathfrak{N}^+}\in\Delta_c\times\Gamma^{-}\cong\mathscr{G}_\infty$ to denote the components of $\sigma_{\mathfrak{N}^+}$ according to the fixed isomorphism $\Delta_c\times\Gamma^{-}\cong\mathscr{G}_\infty$ (note that such $\gamma_{\mathfrak{N}^+}\in\Gamma^-$ is well-defined). We also choose an element $\alpha_c\in\bar{\Q}$ such that $\alpha_c^{-2}=\varphi_t(\sigma_c)\cdot\psi_t(\sigma_c)$. We will also write
\[
\mathcal{R}^-:=\left(\Lambda_\Hf\hat{\otimes}_{\OO_L}\OO_L[\![\Gamma^-]\!]\hat{\otimes}_{\OO_L}\OO_L[\![\Gamma^-]\!]\right)[1/p]
\]
in what follows.    
\end{nota}

\begin{prop}
\label{adjustingfactor}
There exists an element $\mathscr{A}_{\Hf\Hg\Hh}\in\mathcal{R}^-$ such that \begin{itemize}
    \item [(i)] for infinitely many $k\in U_\Hf\cap\Z_{>2}$ and for all $\hat{\nu},\hat{\mu}\in\mathfrak{X}^\mathrm{crit}_{p,k}$, it holds (with $f=\Hf_{\!k}$ as usual)
\[
\mathscr{A}_{\Hf\Hg\Hh}(k,\hat{\nu},\hat{\mu})=\frac{\eta_f}{\lambda_B(k)\cdot\mathcal{E}_p(f,\mathrm{Ad})\cdot\delta_K^{k-1}}\cdot\varphi^-\hat{\nu}(\gamma_{\mathfrak{N}^+})\cdot\psi^-\hat{\mu}(\gamma_{\mathfrak{N}^+})\cdot\frac{\alpha_c}{c\cdot u_K^2}\; ,
\]
\item [(ii)] for all $\hat{\nu},\hat{\mu}\in\mathfrak{X}^\mathrm{crit}_{p,2}$, $
\mathscr{A}_{\Hf\Hg\Hh}(2,\hat{\nu},\hat{\mu})\neq 0$.
\end{itemize}
\begin{proof}
It follows from \cite[lemma 3.3]{BSVast2} that there exists an element $\mathscr{A}_\Hf\in\Lambda_\Hf[1/p]$ such that for infinitely many $k\in U_\Hf\cap\Z_{>2}$ it holds
\[
\mathscr{A}_{\Hf_{\!k}}=\frac{\eta_f}{\lambda_B(k)\cdot\mathcal{E}_p(f,\mathrm{Ad})\cdot\delta_K^{k-1}}\, .
\]
and such that $\mathscr{A}_\Hf(2)\neq 0$. We now set
\[
u:=\frac{\alpha_c\cdot \varphi^-(\gamma_{\mathfrak{N}^+})\cdot \psi^-(\gamma_{\mathfrak{N}^+})}{c\cdot u_K^2}\in L^\times\,.
\]
Then the element $\mathscr{A}_{\Hf\Hg\Hh}:=u\cdot\left(\mathscr{A}_\Hf\hat{\otimes}[\gamma_{\mathfrak{N}^+}]\hat{\otimes}[\gamma_{\mathfrak{N}^+}]\right)\in\mathcal{R}^-$ visibly satisfies the required interpolation property (cf. notation \ref{R-notation}).
\end{proof}
\end{prop}

\begin{defi}
\label{anticycprojdefi}
In the setting \ref{tripleprodass}, the image of $\mathcal{L}^f_p(\Hf,\Hg,\Hh)$ under the map $\mathrm{pr}_{ac}$ of lemma \ref{weighttransitionlemma} is denoted by $\mathcal{L}^f_{p,ac}(\Hf,\Hg,\Hh)$ and called the \emph{anticyclotomic projection} of $\mathcal{L}^f_p(\Hf,\Hg,\Hh)$.
\end{defi}

\begin{teo}
\label{factthm}
Under the natural identification
\[
\mathcal{R}^-=\left(\Lambda_\Hf\hat{\otimes}_{\OO_L}\OO_L[\![\Gamma^-]\!]\hat{\otimes}_{\OO_L}\OO_L[\![\Gamma^-]\!]\right)[1/p]\cong \left(R_{\Hf,\Gamma^-}\hat{\otimes}_{\Lambda_\Hf} R_{\Hf,\Gamma^-}\right)[1/p]\, ,
\]
we have that
\begin{equation}
\label{factoreq}
    \mathcal{L}^f_{p,ac}(\Hf,\Hg,\Hh)=\pm\mathscr{A}_{\Hf\Hg\Hh}\cdot\left(\varphi^-\left(\Theta_\infty^{\mathrm{Heeg}}(\Hf,\varphi_t)\right)\hat{\otimes}\;\psi^-\left(\Theta_\infty^{\mathrm{Heeg}}(\Hf,\psi_t)\right)\right)
\end{equation}
as elements of $\mathcal{R}^-$.
\begin{proof}
It is enough to check that the squares of both sides of equation \ref{factoreq} agree, when specialized to $(k,\hat{\nu},\hat{\mu})$ for infinitely many $k\in U_{\Hf}\cap\Z_>2$ and for every $\hat{\nu}$ and $\hat{\mu}$ finite order characters of $\Gamma^{-}$ (so that $\varphi^-\hat{\nu}$ and $\psi^-\hat{\mu}$ lie in $\mathfrak{X}^\mathrm{crit}_{p,k}$ for every such $k$).

\medskip
We have
\[
\mathcal{L}^f_{p,ac}(\Hf,\Hg,\Hh)(k,\hat{\nu},\hat{\mu})=\mathcal{L}^f_p(\Hf,\Hg,\Hh)(k,y,z)
\]
for any $y,z$ satisfying condition \ref{condyz}.

\medskip
On the other hand we have that
\[
\varphi^-\left(\Theta_\infty^{\mathrm{Heeg}}(\Hf,\varphi_t)\right)(k,\hat{\nu})=\Theta_\infty^{\mathrm{Heeg}}(\Hf,\varphi_t)(k,\varphi^-\hat{\nu})
\]
and
\[
\psi^-\left(\Theta_\infty^{\mathrm{Heeg}}(\Hf,\psi_t)\right)(k,\hat{\mu})=\Theta_\infty^{\mathrm{Heeg}}(\Hf,\psi_t)(k,\psi^-\hat{\mu})\,.
\]

\medskip
The result follows putting together the following ingredients:
\begin{itemize}
    \item [(i)] the factorization of the corresponding complex $L$-functions (cf. equation \eqref{Lfctsfactor1} and lemma \ref{Galoisfact});
    \item [(ii)] the comparison formulas \eqref{interpupgrade} and \eqref{anticycinterp};
    \item [(iii)] our explicit computations for the local factor $\mathscr{I}^{unb}_{\Pi_w,p}$ (cf. proposition \ref{localfactp} and lemma \ref{rootnumberlemma});
    \item [(iv)] the control on the factor $\mathscr{A}_{\Hf\Hg\Hh}$, as described in proposition \ref{adjustingfactor}.
\end{itemize} 
\end{proof}
\end{teo}

\section{Derivatives of triple product \texorpdfstring{$p$}{p}-adic \texorpdfstring{$L$}{L}-functions and Heegner points}
\label{applicationssection}
In this section we describe some applications of theorem \ref{factthm}. We keep the notation as in the previous section (cf. setting \ref{tripleprodass}).

\subsection{Heegner points and Tate's parametrization}
Let $p>3$ denote our fixed prime and let $E/\Q$ be an elliptic curve with multiplicative reduction at $p$. This means that the conductor of $E$ is of the form $N_E=N_E^\circ\cdot p$ with $p\nmid N_E^\circ$. We let $f_E\in S_2(\Gamma_0(N_E))$ to denote the cuspidal newform associated to $E$ via modularity, whose $q$-expansion at $\infty$ will be denoted
\[
f_E=\sum_{n=1}^{+\infty}a_n(E)q^n\; .
\]
In particular we have $a_n(E)\in\Z$ for all $n\geq 1$ and $a_p(E)=1$ (resp. $a_p(E)=-1$) if $E$ has split (resp. non-split) multiplicative reduction at $p$. We write $\alpha:=a_p(E)\in\{\pm 1\}$ in the sequel.

\medskip
Hida theory shows that there exists a unique primitive Hida family
\[
\Hf\in\SSS^{ord}(N_\Hf,\mathbbm{1},\Lambda_\Hf)
\]
of tame level $N_\Hf:=N_E^\circ$ and trivial tame character, such that $\Hf_{\!2}=f_E$.

\medskip
This family will play the role of the Hida family $\Hf$ of the previous section. As for the rest, we keep working in the setting \ref{tripleprodass} and, possibly, add further restrictions. In particular, the conductor $N_E$ of our elliptic curve $E$ is squarefree and satisfies a suitable \emph{Heegner hypothesis} with respect to the fixed quadratic imaginary field $K$.

\medskip
For our applications, we are led to impose one further condition throughout this section.

\begin{ass}
\label{quadraticphi}
$\varphi=\eta_1\eta_2$ has conductor prime to $p$ and $\psi=\eta_1\eta_2^\sigma$ has non-trivial anticyclotomic part $\psi^-$.
\end{ass}

With the notation of section \ref{factorizationsection}, it follows that $\varphi^{-}$ is trivial and that we can identify $\varphi_t=\varphi$.

\medskip
Following the discussion in \cite[section 4.3]{BD2007}, one can define a Heegner point
\begin{equation}
P_\varphi\in \begin{cases} E(H_\varphi)^\varphi & \text{ if }\varphi\neq 1\\
E(K)\otimes\Q & \text{ if }\varphi=1 
\end{cases}
\end{equation}
associated with $\varphi$, essentially coming from a (minimal) parametrisation of $E$ in terms of the Jacobian of a suitable Shimura curve. Here $H_\varphi$ is the field cut out by $\varphi$. Note that, since $p$ is inert in $K$ and $H_\varphi$ is contained in the Hilbert class field of $K$, it follows that $p$ splits completely in $H_\varphi$, so that we can fix an embedding $H_\varphi\subset\Q_{p^2}$ and view the point $P_\varphi$ as a point in $E(\Q_{p^2})\otimes \Q$. Under this identification, the Galois actions on $P_\varphi$ of the Frobenius (as generator of $\Gal(\Q_{p^2}/\Q_p)$) and of any Frobenius element for the abelian extension $H_\varphi/\Q$ coincide. It follows that the points
\[
P^{\pm}_{\varphi,\alpha}:=P_\varphi\pm \alpha\cdot P_\varphi^{\Frob_\p}\in E(H_{\varphi})\otimes\Q\,.
\]
do not depend on the choice of prime $\p$ of $H_\varphi$ above $p$. In what follows, we fix the choice induced by our fixed embedding $\iota_p:\bar{\Q}\hookrightarrow\bar{\Q}_p$ and we view the points $P_\varphi$ and $P^{\pm}_{\varphi,\alpha}$ as elements of $E(\Q_{p^2})\otimes\Q$ under such an embedding.

\medskip
Since $E$ has multiplicative reduction at $p$, it admits a Tate parametrisation, i.e., there is an isomorphism of rigid analyitic varieties
\begin{equation}
\label{Tateunif}
\Phi_\mathrm{Tate}: \mathbb{G}^\mathrm{rig}_{m,\Q_{p^2}}/q_E^{\Z}\xrightarrow{\cong}E_{\Q_{p^2}}^\mathrm{rig}.
\end{equation}

One can define the branch $\log_{q_E}:\C_p^\times\to\C_p$ of the $p$-adic logarithm, uniquely determined by the condition $\log_{q_E}(q_E)=0$, where $q_E\in p\Z_p$ is Tate's $p$-adic period associated with $E$. This yields a logarithm 
\begin{equation}
\log_E:=\log_{q_E}\circ\,\Phi_\mathrm{Tate}^{-1}: E(\Q_{p^2})\to\Q_{p^2}
\end{equation}
at the level of $\Q_{p^2}$-rational points.
\subsection{Restriction to the line \texorpdfstring{$(k,1,1)$}{(k,1,1)}}
We now restrict our attention to the \emph{line} $(k,1,1)$. Recall that $y=1$ (or $z=1$) means that  we consider the specializations given by $y([u])=z([u])=u$ on group-like elements $u\in W_K$. For the first variable, we let $k$ vary in $U_\Hf\cap\Z_{\geq 2}$ (same notation as in remark \ref{anticycchoicesrmk}). The corresponding characters of $\Gamma^-$ via equation \ref{condyz} are clearly both the trivial character $1_{\Gamma^-}$.

\medskip
An easy check shows that, with this choice of specializations, the square of the element 
\[
\mathcal{L}_p(\Hf/K,\varphi):=\Theta_\infty^\mathrm{Heeg}(\Hf,\varphi_t)(\,\cdot\,,1_{\Gamma^-})\in\Lambda
\]
interpolates the algebraic part of the special values $L(\Hf_{\!k}^\circ/K,\varphi,k/2)$, at least when $k>2$. For $k=2$ the $p$-adic multiplier $e_p(f_E,\varphi)$ (cf. proposition \ref{anticyclotomicinterpolationprop}) vanishes, as a manifestation of a so-called \emph{exceptional zero} for our $p$-adic $L$-function. 

Moreover, we see that the element $\mathcal{L}_p(\Hf/K,\varphi)$ coincides with the square-root Hida-Rankin $p$-adic $L$-function attached to $\Hf$ and $\varphi$ in \cite{BD2007}. This follows comparing the above stated interpolation formula \ref{anticycinterp} and the one of  \cite[theorem 3.8]{BD2007}.

\medskip
We can now state one of the main results of \cite{BD2007} (extended to the case of not necessarily quadratic characters $\varphi=\eta_1\eta_2$).

\begin{teo}(\cite[theorem 4.9]{BD2007})
\label{BDHidafamthm}
In the setting described above, it holds
\[
\frac{d}{dk}\mathcal{L}_p(\Hf/K,\varphi)_{|k=2}=\frac{\log_E(P_{\varphi,\alpha}^+)}{2}
\]\,.
\end{teo}

\begin{defi}
We set $\mathcal{L}_p(\Hf/K,\psi):=\psi^-\left(\Theta_\infty^\mathrm{Heeg}(\Hf,\psi_t)\right)(\,\cdot\,,1_{\Gamma^-})\in\Lambda_\Hf$
and we define the restriction to the line $(k,1,1)$ of $\mathcal{L}_p^f(\Hf,\Hg,\Hh)$ as
\[
\mathcal{L}_p^f(\Hf,g,h):=\mathcal{L}_p^f(\Hf,\Hg,\Hh)(\,\cdot\,,1,1)=\mathcal{L}_{p,ac}^f(\Hf,\Hg,\Hh)(\,\cdot\,,1_{\Gamma^-},1_{\Gamma^-})\in\Lambda_\Hf\, .
\]    
\end{defi}

\begin{cor}
\label{derivkcor}
In the above setting (in particular under assumption \ref{quadraticphi}), assume that $L(f_E/K,\psi,1)\neq 0$. Then $\mathcal{L}_p^f(\Hf,g,h)(2)=0$ and
\[
\frac{d}{dk}\mathcal{L}_p^f(\Hf,g,h)_{|k=2}=\frac{c_E}{2}\cdot\log_E(P_{\varphi,\alpha}^+)\,,
\]
where $c_E=\pm\mathscr{A}_{\Hf\Hg\Hh}(2,\hat{\nu}_{1,1},\hat{\mu}_{1,1})\cdot\mathcal{L}_p(\Hf/K,\psi)(2)\in\bar{\Q}_p^\times$. 

In particular, $\frac{d}{dk}\mathcal{L}_p^f(\Hf,g,h)_{|k=2}=0$ if and only if the point $P_{\varphi,\alpha}^+$ is of infinite order.
\begin{proof}
This follows immediately from the above theorem \ref{BDHidafamthm}, the running hypothesis, lemma \ref{adjustingfactor} and the factorization proven in theorem \ref{factthm}. Note that
\[
\mathcal{L}_p(\Hf/K,\psi)(2)=\Theta_\infty^\mathrm{Heeg}(\Hf,\psi_t)(2,\psi^-)\neq 0\,.
\]
Indeed, by assumption \ref{quadraticphi}, we have that $\psi^-$ is non-trivial, so that the $p$-adic multiplier $e_p(\Hf_{\!k},\psi_t\psi^-)$ of the interpolation formula \ref{anticycinterp} never vanishes for $k\in U\cap\Z_{\geq 2}$.
\end{proof}
\end{cor}

\begin{rmk}
Note that (cf. remark \ref{definitenonvanish}) the condition $L(f_E/K,\psi,1)\neq 0$ is \emph{generically} expected to be satisfied.
\end{rmk}

\subsection{Restriction to the line \texorpdfstring{$(2,\nu,\nu)$}{(2,n,n)}}
In this section we fix the weight $k=2$ and we let the anticyclotomic twists vary along the \emph{diagonal} of $\mathfrak{X}^\mathrm{crit}_{p,2}\times\mathfrak{X}^\mathrm{crit}_{p,2}$. In this situation, $\mathfrak{X}^\mathrm{crit}_{p,2}$ is given by finite order characters of $\Gamma^-$.

\begin{defi}
We define the restriction of $\mathcal{L}_p^f(\Hf,\Hg,\Hh)$ to the line $(2,\nu,\nu)$ as
\[
\mathcal{L}_{p,ac}^f(f_E,\Hg\Hh):=\mathcal{L}_{p,ac}^f(\Hf,\Hg,\Hh)_{|k=2,\hat{\nu}=\hat{\mu}}\in\OO_L[\![\Gamma^-]\!]\,.
\]
We also set
\[
\theta_\infty(E/K,\varphi):=\Theta_\infty^\mathrm{Heeg}(\Hf,\varphi_t)_{|k=2}\in\OO_L[\![\Gamma^-]\!]
\]
and
\[
\theta_\infty(E/K,\psi):=\psi^-\left(\Theta_\infty^\mathrm{Heeg}(\Hf,\psi_t)\right)_{|k=2}\in\OO_L[\![\Gamma^-]\!]\, .
\]    
\end{defi}

\medskip
One can check that, under our assumptions, the element $\theta_\infty(E/K,\varphi)$ coincides with the theta-element defined by Bertolini-Darmon (cf. \cite[section 2.7]{BD1996}) in the case of trivial tame character and in more generality by Chida-Hsieh (\cite{CH2018}) and Hung (\cite{Hu2017}). Similarly, the element $\theta_\infty(E/K,\psi)$ is essentially a shift a such a theta-element.

\medskip
Any choice of topological generator $\gamma_0\in\Gamma^-$ gives rise to a topological isomorphism
\begin{equation}
\label{anticycformal}
    \OO_L[\![\Gamma^-]\!]\cong\OO_L[\![T]\!]
\end{equation}
sending $\gamma_0$ to $1+T$. One of the main results of \cite{BD1998} can be stated as follows.
\begin{teo}(cf. \cite[theorem B]{BD1998})
\label{BD98thmB}
The element $\theta_\infty(E/K,\varphi)$ lies in the augmentation ideal of $\OO_L[\![\Gamma^-]\!]$. Equivalently, viewing $\theta_\infty(E/K,\varphi)$ as an element of $\OO_L[\![T]\!]$ via the above identification \ref{anticycformal}, we have
\[
\theta_\infty(E/K,\varphi)\in T\cdot\OO_L[\![T]\!]\,.
\]
Moreover, taking derivatives we obtain
\[
\frac{d}{dT}\theta_\infty(E/K,\varphi)_{|T=0}=\log_E(P_{\varphi,\alpha}^-),
\]
This formula does not depend on the choice of a topological generator of $\Gamma^-$.
\end{teo}

This leads to the following result concerning our triple product $p$-adic $L$-function.
\begin{cor}
\label{derivTcor}
In the above setting (in particular under assumption \ref{quadraticphi}), assume that $L(f_E/K,\psi,1)\neq 0$. View $\mathcal{L}_{p,ac}^f(f_E,\Hg\Hh)$ as an element of $\OO_L[\![T]\!]$ via \eqref{anticycformal}. Then $\mathcal{L}_{p,ac}^f(f_E,\Hg\Hh)_{|T=0}=0$ and
\[
\frac{d}{dT}\mathcal{L}_{p,ac}^f(f_E,\Hg\Hh)_{|T=0}=c_E\cdot\log_E(P_{\varphi,\alpha}^-)\,,
\]
where $c_E\in\bar{\Q}_p^\times$ is the same explicit constant as in corollary \ref{derivkcor}.
\begin{proof}
This follows essentially from the above theorem \ref{BD98thmB}, the factorization of theorem \ref{factthm} and the running hypothesis, in the same way as corollary \ref{derivkcor}.
\end{proof}
\end{cor}

\subsection{A corollary}
Keeping the same setting as in the previous sections (in particular assumption \ref{quadraticphi}), we impose moreover that $\varphi=\varphi_t$ is a quadratic (or genus) character of $K$. 

\medskip
As explained in \cite[section 3.1]{BD2007}, if the quadratic character $\varphi$ is non-trivial, it cuts out a biquadratic extension $H_\varphi=\Q(\sqrt{d_1},\sqrt{d_2})$ where $d_i$ is a fundamental discriminant for $i=1,2$ and $d_1d_2=-d_K$. If we define $\varphi_i$ to be the Dirichlet character attached to the quadratic extension $\Q(\sqrt{d_i})$ for $i=1,2$, one sees that $\varphi_1\varphi_2=\varepsilon_K$. In particular we get
\[
\varphi_1(-N_E)\varphi_2(-N_E)=\varepsilon_K(-N_E)=-1
\]
where the last equality follows from our Heegner assumption.

When $\varphi$ is trivial, one sets $H_\varphi=K$ (this situation corresponds to the case $\{d_1,d_2\}=\{1,-d_K\}$).

\medskip
If $\lambda_E\in\{\pm 1\}$ denotes the eigenvalue relative to $f_E$ for the Atkin-Lehner involution $w_{N_E}$, we can always assume (up to reordering) that
\[
\varphi_1(-N_E)=\lambda_{N_E}\, ,\qquad \varphi_2(-N_E)=-\lambda_{N_E}\,.
\]
Moreover, it follows from \cite[corollary 4.8]{BD2007} that
\begin{equation}
\label{frobactionquadratic}
P_\varphi^{\mathrm{Frob}_\p}=\varphi_1(p)P_\varphi\, .
\end{equation}

Here is a corollary combining the discussion of the previous sections.

\begin{cor}
\label{infiniteorderiff}
In the setting described by assumptions \ref{tripleprodass} and \ref{quadraticphi}, assume that $\varphi=\varphi_t$ is quadratic and that $L(f_E/K,\psi,1)\neq 0$. Then the following facts are equivalent:
\begin{itemize}
    \item [(i)] 
\[
\left(\frac{d}{dk}\mathcal{L}_p^f(\Hf,g,h)_{|k=2}\,,\frac{d}{dT}\mathcal{L}_{p,ac}^f(f_E,\Hg\Hh)_{|T=0}\right)\neq(0,0)
\]
\item [(ii)] The point $P_\varphi$ is of infinite order.
\end{itemize}
\begin{proof}
Equation \ref{frobactionquadratic} shows that, under our assumptions,
\[
P_{\varphi,\alpha}^\pm=\begin{cases}
2\cdot P_\varphi & \text{ if } \varphi_1(p)\alpha=\pm 1\\
0 & \text{ if } \varphi_1(p)\alpha=\mp 1
\end{cases}
\]
Then the result follows immediately from corollaries \ref{derivkcor} and \ref{derivTcor} and the fact that the kernel of $\log_E$ is given by finite order points in $E(\Q_{p^2})$.
\end{proof}
\end{cor}

\printbibliography

\end{document}